\documentclass
{article}
\usepackage{amsmath}
\usepackage{amsfonts}
\usepackage{amssymb}
\usepackage{xspace}

\date{}
\makeatletter \@addtoreset{equation}{section} \makeatother

\newcommand{\adl}{\vspace{1\baselineskip}}

\newenvironment{proof}{\par\noindent{\sc Proof:}
}{\hfill\llap{$\Box$}\vspace{1\baselineskip}\par\noindent}
\newenvironment{proofof}{\par\noindent{\sc Proof}
}{\hfill\llap{$\Box$}\vspace{1\baselineskip}\par\noindent}
\newtheorem{theorem}{Theorem}[section]
\newtheorem{proposition}[theorem]{Proposition}
\newtheorem{lemma}[theorem]{Lemma}
\newtheorem{corollary}[theorem]{Corollary}
\newtheorem{remark}[theorem]{Remark}
\newtheorem{definition}[theorem]{Definition}
\newtheorem{example}[theorem]{Example}
\newcommand{\beq}{\begin{equation}}
\newcommand{\eeq}{\end{equation}}

\newcommand{\ba}{\begin{array}}
\newcommand{\ea}{\end{array}}
\newcommand{\bt}{\begin{theorem}}
\newcommand{\et}{\end{theorem}}
\newcommand{\bp}{\begin{proposition}}
\newcommand{\ep}{\end{proposition}}
\newcommand{\bl}{\begin{lemma}}
\newcommand{\el}{\end{lemma}}
\newcommand{\bc}{\begin{corollary}}
\newcommand{\ec}{\end{corollary}}
\newcommand{\bi}{\begin{itemize}}
\newcommand{\ei}{\end{itemize}}
\newcommand{\ben}{\begin{enumerate}}
\newcommand{\een}{\end{enumerate}}
\newcommand{\bpf}{\begin{proof}}
\newcommand{\epf}{\end{proof}}
\newcommand{\bpff}{\begin{proofof}}
\newcommand{\epff}{\end{proofof}}
\newcommand{\bdf}{\begin{definition}\rm}
\newcommand{\edf}{\end{definition}}
\newcommand{\bex}{\begin{example}\rm}
\newcommand{\eex}{\end{example}}
\def\pri{\hbox to 10pt{\hfil\hbox to 0.4pt{\vrule height5pt width0.4pt
                 depth0pt}\vrule width5pt height0.4pt depth0pt\hfil}}
\newcommand{\TC}{{\rm TC}}

\newcommand{\Div}{{\rm div\,}}
\newcommand{\DIV}{{\rm Div\,}}

\newcommand{\gk}{{\bf{k}}}
\newcommand{\gH}{{\bf H}}
\newcommand{\gK}{{\bf K}}

\newcommand{\calL}{{\mathcal L}}

\newcommand{\EE}{{\mathcal E}}
\newcommand{\gR}{{\mathbb R}}

\newcommand{\Nat}{{\mathbb N}}

\newcommand{\A}{{\mathcal A}}
\newcommand{\D}{{\mathcal D}}
\newcommand{\F}{{\mathcal F}}
\newcommand{\G}{{\mathcal G}}

\newcommand{\M}{{\mathcal M}}
\newcommand{\Ha}{{\mathcal H}}

\newcommand{\BV}{\mathop{\rm BV}\nolimits}

\newcommand{\gnu}{{{\bf{\nu}}_u}}

\def\det{\mathop{\rm det}\nolimits}

\def\Circle{\mathaccent"17}

\newcommand{\wc}{\rightharpoonup}

\newcommand{\Sph}{{\mathbb S}}

\def\mesh{\mathop{\rm mesh}\nolimits}
\newcommand{\ttt}{{\mathfrak t}}
\newcommand{\nnn}{{\mathfrak n}}

\newcommand{\ee}{{\bf e}}
\newcommand{\gn}{{\bf n}}
\newcommand{\gm}{{\bf m}}

\newcommand{\Var}{\mathop{\rm Var}\nolimits}

\let\a=\alpha
\let\be=\beta

\let\d=\delta
\let\e=\varepsilon

\let\vf=\varphi

\let\l=\lambda
\let\m=\mu
\let\n=\nu

\let\p=\pi

\let\r=\rho
\let\s=\sigma
\let\t=\theta
\let\tt=\tau
\let\x=\xi

\let\z=\zeta
\let\vf=\varphi
\let\DD=\Delta
\let\LL=\Lambda
\let\SS=\Sigma

\let\vect=\overrightarrow

\let\wdg=\wedge

\let\wid=\widetilde

\let\pa=\partial
\let\sb=\subset

\let\fa=\forall
\let\tim=\times

\let\sm=\setminus

\let\ul=\underline
\let\ds=\displaystyle

\let\lan=\langle
\let\ran=\rangle

\let\i=\infty
\let\lm=\limits

\title{\Large \bf Bounded variation and relaxed curvature of surfaces}
\author{\it Domenico Mucci and Alberto Saracco
\footnote{%
{\sc Dipartimento di Scienze Matematiche,
Fisiche ed Informatiche, Universit\`{a} di Parma,
Parco Area delle Scienze 53/A, I-43124 Parma, Italy.
E-mail: domenico.mucci@unipr.it, alberto.saracco@unipr.it}
}
}
\begin{document}
\topskip=1.5truecm \maketitle \topskip=1.5truecm \maketitle
%
{\small {\bf Abstract.} We consider a relaxed notion of energy of non-parametric codimension one surfaces that takes into account
area, mean curvature, and Gauss curvature. It is given by the best value obtained by approximation with inscribed polyhedral surfaces.
The BV and measure properties of functions with finite relaxed energy are studied.
Concerning the total mean and Gauss curvature, the classical counterexample by Schwarz-Peano to the definition of area is also analyzed.}
\adl\par\noindent
\adl\par\noindent
\adl\par\noindent
Following the notion of {\em Jordan length} of a curve, the first attempt to define the area of a non-smooth surface \,$\SS$\, was given by J.~A.~Serret in 1868 as
the limit of the elementary area of any sequence of inscribed polyhedral surfaces $P$ converging to \,$\SS$.
The above definition was shown to be incorrect by H.~A.~Schwarz in 1880 (and first published by C.~Hermite in the second edition of his mimeographed lecture notes,
in 1883) and by G.~Peano in 1882 (who published his work in 1890).
In the celebrated example by Schwarz-Peano, they independently showed that if \,$\SS$\, is an ordinary cylinder of radius $R$ and height $H$,
one can define a sequence of inscribed polyhedral surfaces given by the union of congruent triangles with diameters tending to zero, but whose total area converges to any real number not less than the area \,$2\p R H$\, of the cylinder.
\par In the following years, several approaches to provide a correct definition of area were proposed, all based on the principle of lower semicontinuity.
%
The most used is the {\em relaxed area} defined by H.~Lebesgue in 1900. For a codimension one surface \,$\SS$, it is given by the {\em lower limit} of the elementary areas of the polyhedral surfaces uniformly approaching \,$\SS$.
\par In the non-parametric case, the surface \,$\SS$\, is assumed to be the {\em graph}
$$ \G_u=\{(x,u(x))\mid x\in Q\} $$
of a {\em continuous} and real valued function \,$u$\, defined on a closed and bounded domain \,$Q\sb \gR^2$, e.g., \,$Q=[0,1]^2$, the unit square. In his celebrated paper of 1926, L.~Tonelli showed that the graph surface \,$\SS$\, has {\em finite relaxed area}
 in Lebesgue's sense if and only if \,$u$\, is a function of {\em bounded variation}, see \cite{AFP}.
\par The aim of this paper is to extend (at least partially) Tonelli's result concerning the area to a similar notion of {\em total mean and Gauss curvature}.
In correspondence to a relaxed formula that takes into account both area and curvatures, one expects that if \,$u$\, has {\em finite relaxed energy}, then both \,$u$\, and the outward unit normal \,$\n_u$\, are function of bounded variation. Moreover, the non-smooth counterpart of the density of the total mean and Gauss curvature energy of smooth functions \,$u$, suggests that suitable distributions (depending on the approximate derivative of \,$u$, see \cite[Sec.~3.6]{AFP}, and of the unit normal \,$\n_u$\,)
are expected to be measures with {\em finite total variation}, too.
\adl\par In order to tackle the above problem, we recall from J.~M.~Sullivan \cite{Su_sup} the definition of mean curvature and Gauss curvature of a polyhedral surface \,$P$\, in \,$\gR^3$.
\par The mean curvature is supported on the
edges $e$ of $P$, where it is given by $$ \gH_P(e):=\calL(e)\cdot 2\,\sin(\t_e/2) $$
$\calL(e)$\, denoting the {\em length} of the edge and \,$\t_e$\, the exterior {\em dihedral angle} of \,$P$\, along the edge.
\par The Gauss curvature, instead, is supported on the vertexes \,$V$\, of
\,$P$, and in order that the Gauss-Bonnet theorem continues to hold, at each vertex it is given by the {\em angle defect}
$$\gK_{P}(V):=2\pi-\sum_i\t_i $$
where \,$\t_i$\, is the angle of the $i^{th}$-face of \,$P$\, meeting at \,$V$.
Therefore, if e.g. \,$P$\, is the Schwarz-Peano lantern, one has \,$\gK_{P}(V)=0$\, at each vertex, as \,$P$\, is a developable surface.
\par The natural notion of {\em total energy} of \,$P$\, is therefore given by:
$$\EE(P):=A(P)+\EE_\gH(P)+\EE_\gK(P) $$
where \,$A(P)$\, the {\em area} of the polyhedral surface, and \,$\EE_\gH(P)$\, and \,$\EE_\gK(P)$, which will be called the
{\em mean curvature energy} and the {\em Gauss curvature energy}, are respectively defined by
\beq\label{EPint}
 \EE_\gH(P):=\sum_{e\in P} |\gH_P(e)|\,,\qquad  \EE_\gK(P):=\sum_{V\in P} |\gK_{P}(V)| \eeq
where the first summation is taken on all the edges of \,$P$, and the second one
on all the vertexes of \,$P$.
\par We shall consider {\em triangulated} polyhedral surfaces \,$P$\, which are {\em inscribed} in the graph \,$\G_u$\, of a continuous function \,$u:Q\to\gR$,
where \,$Q:=[0,1]^2$\, is the unit square  of \,$\gR^2$.
The relaxed notion of area of the graph of \,$u$\, may be thus written as:
$$ \A(u,Q):=\inf\{\liminf_{h\to\i} A(P_h)\} $$
where {\em the infimum is taken among all the sequences \,$\{P_h\}$\, of inscribed polyhedral surfaces whose corresponding meshes tend to zero}.
Actually, Tonelli's theorem continues to hold: the function \,$u$\, has bounded variation in \,$Q$\, if and only if \,$\A(u,Q)<\i$,
see Proposition~\ref{PTonelli}.
\par In the same spirit, we introduce the {\em relaxed energy} of a continuous function \,$u:Q\to\gR$\, by the formula:
\beq\label{Erelint} \EE(u,Q):=\inf\{\liminf_{h\to\i} \EE(P_h)\} \eeq
where the infimum is taken as above, and the energy \,$\EE(P_h)$\, is given by \eqref{EPint}.
The aim of this paper is to study the $\BV$ and measure properties of the class of functions with finite relaxed energy.
\par We finally point out that a different approach to curvature approximation by polyhedra can be found in \cite{HPW}, where a list of papers on this subject from the point of view of
{\em discrete geometry} is provided.
\adl\par\noindent{\large\sc Outline of the paper.} In Sec.~\ref{Sec:curves}, we collect some features from Sullivan \cite{Su_curv},
concerning the {\em total curvature} of (polygonal) curves.
We shall then prove, Proposition~\ref{Pcurv}, that if a curve has finite total curvature, then the unit normal, when seen as a function of the arc-length parameter,
is a function of bounded variation, with total variation equal to the {\em curvature force}.
\par In Sec.~\ref{Sec:rel}, we introduce our notion of relaxed energy, recalling the definition of mean and Gauss curvature of a polyhedral surface \,$P\sb\gR^3$.
We then see that the Schwarz-Peano counterexample gives a similar drawback concerning the mean curvature:
in general, it does not suffice to take a sequence of polyhedral surfaces inscribed in the cylinder and with diameters of the triangles tending to zero.
Finally, we report the notion by G.~Anzellotti, R.~Serapioni, and I.~Tamanini in \cite{AST} of {\em curvature energy} for smooth surfaces \,$\M$,
and how it is rephrased in the non-parametric case, see \cite{Mu}.
\par In Sec.~\ref{Sec:smoothing}, we analyze the curvature energy of {\em smooth approximations of a polyhedral surface}.
In fact, as it is clear from the converse implication in Tonelli's theorem, in order to obtain the $\BV$-property of a function with finite relaxed area,
one is induced to search for smooth approximating functions.
Concerning the area and the total mean curvature energy \,$\EE_\gH(P)$, a convolution argument yields the expected energy bound for the smooth approximating surfaces,
Proposition~\ref{Pappr}.
However, in general a similar bound of the integral of the modulus of the Gauss curvature of the smooth approximating surfaces cannot be obtained in terms
of the Gauss curvature energy \,$\EE_\gK(P)$. This will be shown in Example~\ref{Ecylinder}, where \,$P$\, is a piece of the Schwarz-Peano lantern.
Roughly speaking, at any vertex \,$V$\, in \,$P$\, we know that \,$\gK_P(V)=0$, whereas in a small neighborhood of each one of the six edges meeting at \,$P$, the outward unit normal of a smooth approximating function has to cover an arc in the Gauss sphere \,$\Sph^2$\, connecting the points given by the values of the outward unit normal to the two triangles of \,$P$\, meeting at the edge.
Therefore, the mapping area of the smooth unit normals gives a contribution equal to the area (with multiplicity) of such a spherical shell in the Gauss sphere \,$\Sph^2$, see also
Remark~\ref{Rvertexes}.
\par In general, a rough area estimate holds, Proposition~\ref{Pestimatearea}. On the other hand, if $V$ is an elliptic vertex of a polyhedral surface $P$,
i.e., if the angle defect at \,$V$\, is positive, we will show that {\em the Gauss curvature can be calculated in terms of a suitable area in the
Gauss sphere}, Proposition~\ref{area=curvature}.
As a consequence, if all the vertexes of the polyhedral graph are of elliptic type, we may
extend Proposition~\ref{Pappr} by obtaining a bound of the integral of the modulus of the Gauss curvature of the smooth approximating surfaces in terms
of the Gauss curvature energy \,$\EE_\gK(P)$, see Corollary~\ref{Capprconv}.
\par In Sec.~\ref{Sec:SPGauss}, we return to the Schwarz-Peano example, showing that the equality given by Proposition~\ref{area=curvature} for elliptic
vertexes of a polyhedral surface drastically fails in this case, where the vertexes are of parabolic type, i.e., with Gauss curvature equal to zero,
see Definition~\ref{Dvert}.
We shall also see that by choosing a different triangulation of $\Sigma$, it turns out that area, mean curvature, and Gauss curvature behave as expected: it suffices to inscribe a prism
\,$Q_n$\, with base a regular $n$-agon in \,$\Sigma$\, and then to triangulate the lateral faces of the prism as we like.
This way, Proposition~\ref{area=curvature} continues to hold.
%
Of course, any approximation procedure has to be done in a smart way, depending on the geometry of the surface \,$\SS$, as in general not all triangulations work properly.
\par In Sec.~\ref{Sec:BV}, we prove, Theorem~\ref{TBV}, that if \,$u$\, is a continuous function with finite relaxed energy \eqref{Erelint}, then {\em the outward unit normal \,$\n_u:Q\to\Sph^2$\, is a function of bounded variation}. We remark that the unit normal is well defined a.e. on \,$Q$\, in terms of the
approximate partial derivatives of \,$u$, as \,$u$\, is a function in \,$\BV(Q)$, by Proposition~\ref{PTonelli}.
\par As the case of graphs of smooth functions suggests, an extra term should be added in order to bound the (relaxed) energy corresponding to the mean curvature.
For this purpose, we recall that the {\em distributional divergence} of an $L^1$-vector field \,$\s:Q\to\gR^2$\,
is well-defined by duality through the formula
$$ \lan\DIV \s,\vf\ran:= -\int_Q\s(x)\bullet \nabla\vf(x)\,dx\,,\qquad \vf\in C^\i_c(\Circle Q)$$
where \,$\bullet$\, denotes the scalar product in \,$\gR^2$. If e.g. \,$u:Q\to\gR$\, is a continuous function with finite relaxed energy,
the vector fields \,$\s^j_u:Q\to\gR^2$
$$ \s^j_u:=(-\n_u^j\,\pa_2 u,\,\n_u^j\,\pa_1 u)\,,\quad j=1,2,3 $$
are summable, see Example~\ref{Esju}.
\par When \,$u$\, is smooth, say of class \,$C^2$, it turns out that the distribution \,$\DIV\s^j_u$\, is an absolute continuous signed measure
with density equal to the pointwise divergence of \,$\s^j_u$. Moreover, we have \,$\Div\s^j_u(x)=\m^j_u(x)$\, for each \,$x\in Q$,
where
\beq\label{detint}
\m^j_u(x)=\det\left(\ba{cc}
  \pa_1u(x) & \pa_2u(x) \\
  \pa_1\n_{u}^j(x) & \pa_2\n_{u}^j(x)
\ea\right)\,. \eeq
\par In case of polyhedral surfaces, we in fact see, Proposition~\ref{Pmeanpoly}, that the energy term \,$\EE_\gH(P)$\, in \eqref{EPint} can be seen as the total variation of the vector-valued measure \,$\gm_u:=(\gm_u^1,\gm_u^2,\gm_u^3)$, where
$$ \gm_u^j:= (D\n_u^j,\DIV\s^j_u)\,,\qquad j=1,2,3\,. $$
\par More generally, we prove, Theorem~\ref{Tmeas}, that if a continuous function \,$u$\, has finite relaxed energy \eqref{Erelint}, then
{\em the distributional divergence \,$\DIV \s^j_u$\, is a finite measure}, that is decomposed as
$$
 \DIV \s^j_u=\m^j_u\,\calL^2\pri Q+(\DIV \s^j_u)^s
$$
where \,$\calL^2$\, is the Lebesgue measure in $\gR^2$, the summable function \,$\m^j_u(x)$\, is defined for \,$\calL^2$-a.e. $x\in Q$\, by \eqref{detint}, and \,$(\DIV \s^j_u)^s$\, is singular w.r.t. the Lebesgue measure.
\par As to the Gauss curvature energy of polyhedral surfaces, we do not have an analogous to Proposition~\ref{Pmeanpoly}.
However, as a consequence of Corollary~\ref{Capprconv} we shall obtain, Proposition~\ref{PGauss}, that
if \,$u$\, is a strictly convex function with finite relaxed energy,
then {\em all the $2\tim 2$-minors of the matrix}
\beq\label{matint}  \left(
      \begin{array}{cc}
        \pa_1\n_u^1 & \pa_2\n_u^1 \\
        \pa_1\n_u^2 & \pa_2\n_u^2 \\
        \pa_1\n_u^3 & \pa_2\n_u^3 \\
      \end{array}
    \right)
\eeq
of the approximate partial derivatives of the unit normal
{\em are summable functions}.
\adl\par\noindent{\large\sc Open questions.}
We expect the claim in Proposition~\ref{PGauss} to hold true without assuming strict convexity. However, we are not able to prove this fact, due to the drawbacks illustrated in Example~\ref{EGauss}.
\par On the other hand, it is an open problem to {\em characterize the class \,$\EE(Q)$\, of continuous functions \,$u:Q\to\gR$\, with finite relaxed energy
\eqref{Erelint}}. Starting from our results, one may conjecture that \,$u\in\EE(Q)$\, if and only if the following properties hold:
\ben \item $u$ is a function in \,$\BV(Q)$;
\item the outward unit normal \,$\n_u$\, is a function of bounded variation;
\item for \,$j=1,2,3$, the distributions \,$\DIV\s^j_u$\, are measures with finite total variation;
\item the $2\tim 2$-minors of the matrix \eqref{matint} are summable functions in $L^1(Q)$.
\een

\section{BV-property of a curve with finite total curvature}\label{Sec:curves}
In this section we collect some notions and properties from Sullivan \cite{Su_curv}, concerning the total curvature of (polygonal) curves in Euclidean spaces.
We then prove, Proposition~\ref{Pcurv}, that if a curve has finite total curvature, then the unit normal, when seen as a function of the arc-length parameter,
is a function of bounded variation, with total variation equal to the curvature force.
\par Even if the following statements hold true in high codimension, for our purposes we restrict to consider curves \,$c$\, in \,$\gR^2$\, parameterized by
\,$c:I\to\gR^2$, where \,$I:=[0,1]$\, and
\,$c(t)=(c^1(t),c^2(t))$\, is continuous.
%
\adl\par\noindent{\large\sc Length.} Any polygonal curve \,$P$\,
inscribed in \,$c$, say \,$P \ll c$, is obtained by choosing a
finite partition \,$\D:=\{0=t_0<t_1<\ldots<t_{n-1}<t_{n}=1\}$\, of
\,$I$, say \,$P=P(\D)$, and letting \,$P:I\to\gR^2$\, such that
\,$P(t_i)=c(t_i)$\, for \,$i=0,\ldots,n$, and \,$P(t)$\, affine on
each interval \,$I_i:=[t_{i-1},t_i]$\, of the partition. Setting
\,$\ee_i=\dot P(t)\in\gR^2$\, for \,$t\in \Circle I_i$\, we have
\,$\ee_i\neq 0_{\gR^2}$\, for each \,$i=1,\ldots,n$\, and hence
the length of \,$P$\, is
$$ \calL(P)=\sum_{i=1}^{n}\calL(I_i)\cdot |\ee_i|=\int_I|\dot P(t)|\,dt\,. $$
The length \,$\calL(c)$\, of \,$c$\, is defined by
\,$\calL(c):=\sup\{\calL(P)\mid P\ll c\}$, and \,$c$\, is said to
be {\em rectifiable} if \,$\calL(c)<\i$.
With the above notation, we let
$$\mesh \D:=\sup_{1\leq i\leq n}\calL(I_i)\,,\quad \mesh P:=\sup_{1\leq i\leq n}\calL(I_i)\cdot |\ee_i| \,. $$
By uniform continuity of
\,$c\in C^0(I,\gR^2)$, for each \,$\e>0$\, we can find \,$\d>0$\,
such that \,$\mesh P<\e$\, if \,$\mesh \D<\d$\, and \,$P=P(\D)$. As a
consequence, taking a sequence \,$P_h=P(\D_h)$\, where
\,$\{\D_h\}$\, is any sequence of partitions of \,$I$\, such that
\,$\mesh \D_h\to 0$,
we get \,$\mesh P_h\to 0$\, and hence the
convergence \,$\calL(P_h)\to\calL(c)$\, of the length functional.
\adl\par\noindent{\large\sc Total variation.} Following e.g. \cite[Sec.~3.2]{AFP}, given a (not
necessarily continuous) function \,$f:I\to\gR^2$\, with finite pointwise variation, the (essential) {\em total variation}
\,$\Var_{\gR^2}(f)$\, is the infimum of the pointwise variation computed among the functions that agree with $f$ at $\calL^1$-a.e. $t\in I$, where $\calL^1$ is the Lebesgue measure in $\gR$. The function $f$ is bounded and summable in $I$,
and its distributional derivative \,$Df$, given by
$$ \lan Df,\vf\ran:=-\int_I f(t)\bullet \dot\vf(t)\,dt\,,\qquad
\vf\in C^\i_c(\Circle I,\gR^2)\,.$$
is a finite measure. Moreover, one has \,$\Var_{\gR^2}(f)=|Df|(I)$, where
$$ |Df|(I):=\sup\{\lan Df,\vf\ran\mid \vf\in C^\i_c(\Circle I,\gR^2)\,,\,\,\Vert\vf\Vert_\i\leq 1\}<\i\,.$$
If \,$f\in\BV(I,\gR^2)$, the approximate derivative \,$\dot f$\, is an \,$L^1$-function, and one may decompose the distributional derivative into its (mutually singular)
absolutely continuous, jump, and Cantor components, respectively:
$$ Df=D^af+D^Jf+D^Cf $$
where the absolutely continuous component reads as \,$D^af=\dot f\,\calL^1\pri I$,
the Jump component \,$D^Jf$\, is concentrated on an at most countable subset of \,$I$,
and the Cantor component is a diffuse part, so that \,$D^Cf(A)=0$\, if \,$\Ha^0(A)<\i$.
In particular, a (continuous) curve \,$c$\, as above
is rectifiable if and only if \,$c\in\BV(I,\gR^2)$, and
in this case \,$\calL(c)=\Var_{\gR^2}(c)=|Dc|(I)$.
\par Assume now that \,$|f(t)|=1$\, for $\calL^1$-a.e. \,$t\in
I$, i.e., \,$f$\, is a measurable function from \,$I$\, into the unit circle
\,$\Sph^1:=\{y\in\gR^2\,:|y|=1\}$\, of \,$\gR^2$. The (essential) total
variation of \,$f$\, can be computed in two different ways, by
taking the geodesic distance in \,$\Sph^1$\, or the Euclidean
distance in \,$\gR^2$. Since \,$d_{\gR^2}(Q_1,Q_2)\leq
d_{\Sph^1}(Q_1,Q_2)\leq (\p/2)\cdot d_{\gR^2}(Q_1,Q_2)$\, for any
\,$Q_1,Q_2\in\Sph^1$, in general one obtains:
$$ \frac{2}{\p}\,\Var_{\Sph^1}(f)\leq \Var_{\gR^2}(f)\leq \Var_{\Sph^1}(f) $$
and hence {\em \,$f$\, has bounded total variation in \,$\Sph^1$\,
if and only if it has bounded total variation in \,$\gR^2$}, i.e.
$$ f\in \BV_{\Sph^1}\iff f\in \BV_{\gR^2}\,.$$
In this case, with a modern notation one writes \,$f\in\BV(I,\Sph^1)$.
In particular, if \,$f$\, is smooth one has
\,$\Var_{\gR^2}(f)=\Var_{\Sph^1}(f)=\int_I|\dot f(t)|\,dt$, whereas in general the
strict inequality \,$\Var_{\gR^2}(f)<\Var_{\Sph^1}(f)$\, holds, as
a gap appears at each jump point of \,$f$, and we recall that
\,$|Df|(J)=\Var_{\gR^2}(f)$.
\begin{remark} In Example~\ref{ECantor} below, where \,$f$\, is the unit normal to the Cartesian curve given by the graph of a primitive of the Cantor-Vitali function, it
turns out that the
Cantor component \,$D^Cf$\, of the distributional derivative of
\,$f$\, does not produce a gap between the two definitions of
total variation. For this reason we expect that for a function \,$f\in \BV(I,\Sph^1)$, one has
\,$\Var_{\gR^2}(f)=\Var_{\Sph^1}(f)$\, if and only if the Jump
component \,$D^Jf=0$, i.e., if and only if \,$f$\, has a continuous representative. \end{remark}
\adl\par\noindent{\large\sc Total curvature.} Following Milnor
\cite{Mi}, the {\em total curvature} of a curve \,$c$\, is given by
$$\TC(c):=\sup\{\TC(P)\mid P \ll c\} $$
where the total curvature \,$\TC(P)$\, of the inscribed polygonal
\,$P$\, is the sum of the {\em turning angles} \,$\t_i$\, at the
edges of \,$P$. Therefore, denoting by \,$\bullet$\, the scalar
product in \,$\gR^2$, with the above notation we get
$$ \TC(P)=\sum_{i=1}^{n-1} \t_i\,,\qquad
\t_i:=\arccos{\ee_i\bullet \ee_{i+1} \over |\ee_i|\cdot
|\ee_{i+1}| }\,,\quad i=1,\ldots,n-1 $$
(where a further turning angle between $\ee_n$ and $\ee_1$ appears if $P$ is closed) and hence \,$\TC(P)$\, agrees with the total variation
in \,$\Sph^1$\, of the {\em tantrix} (or tangent indicatrix)
\,$\ttt_P$\, (the tantrix assigns to a.e. point the oriented unit
tangent vector in \,$\Sph^1$), i.e.
$$ \TC(P)=\Var_{\Sph^1}(\ttt_P)\,.$$
\par If \,$c$\, has finite total curvature, \,$\TC(c)<\i$, then
\,$c$\, is rectifiable, hence its arc-length parameterization is Lipschitz-continuous. Therefore, by Rademacher theorem \cite[Sec.~2.3]{AFP}, the tantrix
\,$\ttt_c$\, is well-defined a.e. by the derivative of \,$c$\, with respect to its arc-length parameter. Moreover, the total curvature agrees with the (essential) total variation in \,$\Sph^1$ of the tantrix. Therefore, letting
\,$\nnn_c:=\ttt_c^\perp$, where we have set
\,$(a,b)^\perp:=(b,-a)$, we get
$$ \TC(c)=\Var_{\Sph^1}(\ttt_c)=\Var_{\Sph^1}(\nnn_c)\,. $$
Also, taking any sequence \,$P_h=P(\D_h)\ll c$\, with
\,$\mesh\D_h\to 0$, we get convergence \,$\TC(P_h)\to\TC(c)$\, of the total curvature functional.
\adl\par\noindent {\large\sc The curvature force.} The {\em
curvature force} \,$\TC^*(P)$\, of a polygonal is given by the
total variation in \,$\gR^2$\, of the tantrix \,$\ttt_P$\,:
$$ \TC^*(P):=\Var_{\gR^2}(\ttt_P) $$
compare \cite{Su_curv}.
In particular, if \,$P\ll c$, with the previous notation one has:
$$ \TC^*(P)=\sum_{i=1}^{n-1} 2\sin({\t_i/2}) $$
and therefore
$$ {2\over\pi}\,\TC(P)\leq \TC^*(P)\leq \TC(P)\,. $$
Furthermore, we have \,$\dot P\in L^1(I,\gR^2)$\, and the unit normal is
well defined outside the edges of \,$P$\, by
$$ \gn_P(t):={\dot P(t)^\perp \over \vert \dot P(t)\vert}\,,\qquad
t\neq t_i \,. $$
Then \,$\gn_P$\, is a (piecewise constant) function of bounded
variation in the class \,$\BV(I,\Sph^1)$,
and furthermore
$$ |D\gn_P|(I)=\TC^*(P)<\i\,. $$
\par Defining by the same approach as above the Euclidean total curvature, or {\em
curvature force}, of \,$c$\, by
$$\TC^*(c):=\sup\{\TC^*(P)\mid P \ll c\} $$
one infers that \,$c$\, has finite curvature force if and only if
it has finite total curvature. In this case, moreover, taking
again any sequence \,$P_h=P(\D_h)\ll c$\, with \,$\mesh\D_h\to 0$,
one gets the convergence \,$\TC^*(P_h)\to\TC^*(c)$\, of the curvature forces.
\par Finally, if \,$c$\, is smooth, say \,$c\in C^2(I,\gR^2)$, then one has
$$ \calL(c)=\int_I |\dot c(t)|\,dt\,,\qquad
\TC(c)=\TC^*(c)=\int_I{|\dot c\wedge \ddot c| \over |\dot
c|^2}\,dt
$$
where \,$|\dot c\wedge \ddot c|=|\dot c^1\,\ddot c^2-\dot
c^2\,\ddot c^1|$\, if \,$c=(c^1,c^2)$. In fact, denoting by
\,$\gn_c(t)$\, the unit normal at \,$c(t)$\, one gets:
$$ \gn_c(t)={\dot c(t)^\perp \over \vert \dot
c(t)\vert}\,,\qquad |\dot\gn_c(t)|={|\dot c\wedge \ddot c| \over
|\dot c|^2}(t)\qquad  \fa\, t\in I\,. $$
\par\noindent{\large\sc $BV$-property.} Let \,$c$\, be a rectifiable curve, so
that \,$L:=\calL(c)<\i$. Let \,$\ul c:I_L\to\gR^2$\, be the
arc-length parameterization of \,$c$, where \,$I_L:=[0,L]$, so
that \,$\dot{\ul c}\in L^\i(I_L,\gR^2)$\, with \,$|\dot{\ul c}(s)|=1$\,
for a.e. \,$s\in I_L$. Define
$$  \ul\gn_c(s):= { \dot{\ul c}(s)^\perp \over \vert \dot
{\ul c}(s)\vert}=\dot{\ul c}(s)^\perp\,,\qquad s\in I_L\,.
$$
In the following result we recover the definition
\,$\TC^*(c):=\Var_{\Sph^1}(\ttt_c)$\, by Sullivan \cite{Su_curv},
exploiting the $\BV$-property of the unit normal \,$\ul\gn_c$.
\bp\label{Pcurv} If \,$\TC^*(c)<\i$, then \,$\ul\gn_c$\, is a function of
bounded variation in \,$\BV(I_L,\Sph^1)$, and its total variation in $\gR^2$
is equal to the curvature force
and to the total variation in $\gR^2$ of the tantrix \,$\ttt_c$,
i.e.
$$ |D\ul\gn_c|(I_L)=\TC^*(c)=\Var_{\gR^2}(\ttt_c)=\Var_{\gR^2}(\nnn_c)\,. $$ \ep
\bpf Choose \,$(P_h)$\, to be a sequence of polygonals inscribed in
\,$c$\, such that \,$\mesh P_h \to 0$, so that both
\,$L_h:=\calL(P_h)\to L:=\calL(c)$\, and \,$\TC^*(P_h)\to
\TC^*(c)$. Let \,$\vf_h:I_L\to I$\, be the inverse of the
bijective and increasing function \,$\psi_h:I\to I_L$\,
$$\psi_h(t):= {L\over L_h}\int_0^t |\dot P_h(\l)|\,d\l\,,\qquad
t\in I\,. $$
Letting \,$c_h(s):=P_h(\vf_h(s))$, \,$s\in I_L$, we have \,$|\dot
c_h(s)|\equiv L_h/L$\, a.e.,  and hence (by Ascoli-Arzela's
theorem) possibly passing to a subsequence \,$c_h$\, uniformly
converges to some function \,$f\in C^0(I_L,\gR^2)$.
We have \,$f=\ul c$. In fact, using that \,$\mesh P_h\to 0$\, and
\,$\calL(c)<\i$, we deduce that
\,$\psi_h(t)\to\calL(c_{\vert[0,t]})$\, as \,$h\to\i$\, for each
\,$t\in I$. By Dini's theorem we get uniform convergence of
\,$\{\psi_h\}$\, on \,$I_L$\, and on the other hand
$$ c_h(s)=P_h(t) \iff  s=\psi_h(t)\,,\qquad \ul c(s)=c( \calL(c_{\vert[0,t]})) =\lim_{h\to\i}c(\psi_h(t))\,.$$
%
%
\par Setting now $$\gn_{h}(s):={\dot c_h(s)^\perp\over |\dot
c_h(s)|}={L\over L_h}\,\dot c_h(s)^\perp $$
by definition of curvature force we have
\,$|D\gn_h|(I_L)=\TC^*(P_h)$, whence \,$|D\dot
c_h|(I_L)=(L_h/L)\,\TC^*(P_h)$, with \,$(L_h/L)\to 1$\, and
\,$\TC^*(P_h)\to\TC^*(c)<\i$.
Therefore, we deduce that a subsequence of \,$\{\dot c_h\}$\,
weakly converges in the $\BV$-sense to some function
\,$v\in\BV(I_L,\gR^2)$.
\par We claim that \,$v=\dot{\ul c}$\, a.e. in \,$I_L$, which
clearly yields that the whole sequence \,$\{\dot c_{h}\}$\, weakly
converges to \,$\dot{\ul c}$. In fact, using that by
Lipschitz-continuity
$$ c_h(s)=c_h(0)+\int_0^s \dot c_h(\l)\,d\l\qquad \fa\, s\in I_L
$$
where \,$c_h(0)=c(0)=\ul c(0)$\, for each \,$h$, and setting
$$ V(s):=\ul c(0)+\int_0^s v(\l)\,d\l\qquad \fa\, s\in I_L $$
by the weak $\BV$ convergence \,$\dot c_h\wc v$, which implies the
strong \,$L^1$-con\-ver\-gen\-ce, we have \,$c_h\to V$\, in \,$L^\i$,
hence \,$\dot c_h\wc \dot V=v$\, a.e. in \,$I_L$. But we already
know that \,$c_h\to \ul c$\, in \,$L^\i$, thus \,$v=\dot{\ul c}$.
\par The weak $\BV$ convergence \,$\dot c_h\wc \dot{\ul c}$\, clearly implies
the weak $\BV$ convergence \,$\gn_{h}\wc \ul\gn_c$, whence
\,$\ul\gn_c\in \BV(I_L,\gR^2)$. Then by lower semicontinuity
$$ |D\ul\gn_c|(I_L)\leq\liminf_{h\to\i}|D\gn_h|(I_L)=\lim_{h\to\i}\TC^*(P_h)=\TC^*(c)\,.$$
Finally, arguing as in Sullivan \cite{Su_curv}, one obtains
$$ \lim_{h\to \i}|D\gn_{h}|(I_L)=\lim_{h\to \i}\Var_{\gR^2}(\nnn_{P_h})=\Var_{\gR^2}(\nnn_{c})= |D\ul\gn_c|(I_L) $$
and hence \,$|D\ul\gn_c|(I_L)=\TC^*(c)$. \epf
\bex\label{ECantor} Let \,$c(t)=(t,u(t))$\, be the Cartesian curve
given by the graph of the primitive \,$u(t):=\int_0^tv(\l)\,d\l$\,
of the classical Cantor-Vitali function \,$v:[0,1]\to \gR$\,
associated to the ``middle thirds" Cantor set. We have
\,$\calL(c)=\int_I\sqrt{1+v^2(t)}\,dt<\i$\, and
$$\gn_c(t):={\dot c(t)^\perp\over |\dot c(t)| } = {(v(t),-1) \over \sqrt{1+v^2(t)} } $$
so that
$$ D\gn_c={(1,v)\over{(1+v^2)^{3/2}}}\,D^Cv\,,$$
$$ |D\gn_c|(I)=\int_I{1\over\sqrt{1+v^2(t)}}\,d|D^Cv|=|D\arctan v|(I)\,. $$
\par We now choose the polygonal \,$P_k\ll c$\, corresponding to the
subdivision
$$D_k:=\{t^k_h=h\,3^{-k}\mid h=0,\ldots,3^k\}\,, \qquad k\in\Nat\,. $$
The corner points of \,$P_k$\, agree with the values
\,$(t^k_h,u_k(t^k_h))$, where \,$u_k(t):=\int_0^t v_k(\l)\,d\l$\,
and \,$v_k$\, is the classical $k$-th approximation of the
Cantor-Vitali function \,$v$. Therefore, the polygonal \,$P_k$\, contains \,$2^{k+1}-2$\, vertexes, each edge has slope greater than \,$6^{-k}$, and
the difference between the slopes of two consecutive edges is smaller than \,$2^{-(k+1)}$. As a consequence, each turning angle of \,$P_k$\, is
smaller than \,$\arcsin(4^{-k})$. Using that \,$0\leq\t-2\arcsin(\t/2)\leq (\a/2)^3$\, if \,$0<\t<\a$, we thus get
$$ \lim_{k\to\i}|\TC^*(P_k)-\TC(P_k)|\leq \lim_{k\to\i} 2^{k}\cdot \arcsin(4^{-k})^3=0\,. $$
As a consequence, we obtain
$$ \Var_{\Sph^1}(\gn_c)=\Var_{\gR^2}(\gn_c)=|D\gn_c|(I) $$
and also, compare \cite{AcMu}, that
$$
\TC^*(c)=\TC(c)=|D\gn_c|(I)=\int_I{1\over\sqrt{1+v^2(t)}}\,d|D^Cv|\,.
$$
\par In conclusion, in this example one sees that the occurrence of a Cantor-part in the
derivative of the unit normal \,$\gn_c$\, does not change the
computation when considering the total variation in \,$\Sph^1$\,
or in \,$\gR^2$. \eex
\section{The relaxed energy}\label{Sec:rel}
In this section we introduce the notion of relaxed energy of a continuous function \,$u$, that takes into account the area, the mean curvature, and of the Gauss curvature of the triangulated polyhedral surfaces inscribed in the graph-surface of \,$u$. For this purpose,
we first recall from \cite{Su_sup} the notion of mean and Gauss curvature of a polyhedral surface \,$P\sb\gR^3$.
\par We then see that the classical Schwarz-Peano counterexample
to the definition of area given by Serret, gives a similar drawback concerning the mean curvature:
in order to have a good definition, similarly to what happens for the area, it does not suffice to take a sequence of polyhedral surfaces inscribed in the cylinder and with diameters of the triangles tending to zero.
For this reason, we propose a relaxed definition in the same spirit as in Lebesgue's definition of area.
\par Finally, we recall the notion by Anzellotti-Serapioni-Tamanini in \cite{AST} of {\em curvature energy} for smooth surfaces \,$\M$,
and how it is rephrased in the non-parametric case, i.e. when \,$\M$\, is given by the graph \,$\G_u$\, of a smooth function
 \,$u:Q\to\gR$, see \cite{Mu}.
In the sequel we shall thus consider functions \,$u$\, defined on the unit square \,$Q:=[0,1]^2$\, of \,$\gR^2$.
\adl\par\noindent{\large\sc Polyhedral surfaces.} The {\em mean curvature of a polyhedral surface} \,$P$\, in \,$\gR^3$\, was defined
by Sullivan \cite{Su_sup} in such a way that it is supported on the
edges. Namely, if \,$e$\, is
an edge of \,$P$, then
$$ \gH_P(e):=\calL(e)\cdot 2\,\sin(\t_e/2) $$
where \,$\calL(e)$\, is the length of the edge and \,$\t_e$\, is the exterior dihedral angle of \,$P$\, along the edge.
\par The {\em Gauss curvature of a polyhedral
surface} was defined by Sullivan \cite{Su_sup} in such a way that
the Gauss-Bonnet theorem continues to hold. It is concentrated at
the vertexes, and in the case of a triangulated polyhedral surface
\,$P$, the Gauss curvature at a vertex \,$V$\, agrees with the
{\em angle defect}, whence
$$\gK_{P}(V):=2\pi-\sum_i\t_i $$
where \,$\t_i$\, is the angle of the $i^{th}$-triangle of \,$P$\,
meeting at \,$V$.
\par We then respectively define the {\em mean curvature energy} and the {\em Gauss curvature energy} of a polyhedral surface \,$P$\, by
\beq\label{EPHK}
 \EE_\gH(P):=\sum_{e\in P} |\gH_P(e)|\,,\qquad  \EE_\gK(P):=\sum_{V\in P} |\gK_{P}(V)| \eeq
where the first summation is taken on all the edges of \,$P$, and the second one
on all the vertexes of \,$P$. Denoting by \,$A(P)$\, the {\em area} of the polyhedral surface, we define the {\em the total energy} of \,$P$\, by:
\beq\label{EP} \EE(P):=A(P)+\EE_\gH(P)+\EE_\gK(P)\,. \eeq
\par Finally, for future use, with the above notation we also denote
\beq\label{EPH2}
 \wid\EE_\gH(P):=\sum_{e\in P} \calL(e)\cdot\t_e\eeq
so that we clearly have \,$ \EE_\gH(P) \leq \wid\EE_\gH(P)\leq (\p/2)\cdot \EE_\gH(P)$.
\adl\par\noindent
{\large\sc Relaxed energy.} Let \,$u:Q\to\gR$\, be a continuous function. 
We say that a polyhedral surface is inscribed in the graph of \,$u$\, if we can find a finite triangulation \,$\D$\, of the square domain \,$Q$\, such that \,$P$\, is the graph of the piecewise affine and
continuous function \,$v=v(\D)$\, that agrees with \,$u$\, on the $0$-skeleton of the triangulation and \,$v$\, is affine on each triangle \,$\DD$\, of
the triangulation. In this case we shall write \,$P=P(u,\D)$\, to outline the dependence of \,$P$\, on the values of \,$u$\, on the $0$-skeleton of the
triangulation \,$\D$. The mesh of the triangulation, say \,${\text{mesh}\,}\D$, is given by the supremum of the diameter of the triangles \,$\DD$\, of \,$\D$.
\par We introduce the following relaxed notion of area of the graph of \,$u$\,:
\beq\label{Area} \A(u,Q):=\inf\{\liminf_{h\to\i} A(P_h)\} \eeq
where {\em the infimum is taken among all the sequences \,$\{P_h\}$\, of inscribed polyhedral surfaces whose corresponding meshes tend to zero}; i.e., if \,$P_h=P(u,\D_h)$, then \,${\text{mesh}\,}\D_h\to 0$\, as \,$h\to\i$. By uniform continuity, in fact, this condition implies
that the sequence \,$v_h=v(\D_h)$\, converges to \,$u$\, uniformly on \,$Q$.
\par In the same spirit, we introduce the {\em relaxed energy} of a continuous function \,$u:Q\to\gR$\, by the formula:
\beq\label{Erel} \EE(u,Q):=\inf\{\liminf_{h\to\i} \EE(P_h)\} \eeq
where the infimum is taken as in formula \eqref{Area} above, and \,$\EE(P_h)$\, is given by \eqref{EP}.
\adl\par\noindent{\large\sc Schwarz-Peano example.}\label{SchPea} Consider the lateral surface \,$\SS$\, of a cylinder of radius $R$ and height $H$. Its area is $2\pi R H$, the principal curvatures are $\gk_1=0$ and $\gk_2=1/R$,
whence the Gauss curvature is zero and the integral of the mean curvature \,$\gH=1/(2R)$\, is \,$\int_\SS \gH\,d\Ha^2=\p\,H$.
\par In the classical Schwarz-Peano counterexample to the definition of area given by Serret, one considers for each \,$m,n\in\Nat^+$\, the polyhedral surface \,$P_{m,n}$\, inscribed in \,$\SS$\, and given by \,$2mn$\,
congruent isosceles triangles. With the parameterization  \,$[0, 2\p]\tim[0,H]\ni(\t,z)\mapsto(R\cos\t,R\sin\t,z)$, when \,$m$\, and \,$n$\, are even, the triangles are obtained by taking the
vertexes at the points corresponding to \,$(\p i/m,(H/n)j)$, when both $i$ and $j$ are even, or when both $i$ and $j$ are odd, $i=0,\ldots,m$, $j=0,\ldots,n$.
Letting \,$\a_m:=\p/m$, each triangle has basis \,$b:=2R\sin\a_m$\, and height \,$h:=((H/n)^2+d^2)^{1/2}$, where \,$d:=R(1-\cos\a_m)=2R\sin^2(\a_m/2)$.
Therefore, the area of the polyhedral surface is
$$ A(P_{m,n})=2mnR\sin(\a_m)\,\sqrt{(H/n)^2+4R^2\sin^4(\a_m/2)}\,. $$
We have \,$A(P_{m,n})\geq 2mR\sin(\a_m)\,H$, which tends to \,$2\p RH$\, as \,$m\to \i$. 
Furthermore, when e.g. \,$n=m^2$\, one has \,$A(P_{m,m^2})\to 2\p R\sqrt{H^2+R^2\p^4/4}$, and when \,$n=m^4$\, one gets \,$A(P_{m,m^4})\to +\i$.
However, when \,$m=n^p$\, for some positive integer exponent \,$p$, one obtains convergence to the area of the cylinder, i.e., \,$A(P_{n^p,n})\to 2\p RH$\, as \,$n\to\i$.
\par We now wish to give a similar computation concerning the mean curvature.
As to the Gauss curvature, in fact, we observe that at each interior vertex of \,$P_{m,n}$\, six triangles meet, four ones with an angle \,$\a$\, and two ones with an angle $2\be$, where \,$\a+\be=\p/2$, whence the Gauss curvature at each vertex is \,$2\p-(4\a+2\cdot 2\be)=0$, the polyhedral surface being developable, too.
\par As to the mean curvature of \,$P_{m,n}$, we recall that it is concentrated at the edges \,$e$\, of the triangles, and
at each edge the contribution is given by \,$\calL(e)\cdot\t_e/2$, where \,$\calL(e)$\, is the length and \,$\t_e$\, is the dihedral angle of the two faces meeting at the edge \,$e$, see also Remark~\ref{Rcurv} below. Notice that the factor \,$1/2$\, is due since we define the mean curvature of a smooth surface as \,$\gH=(\gk_1+\gk_2)/2$.
We shall prove the following:
\bp\label{Pex} The total mean curvature of \,$P_{m,n}$\, converges to the
integral on \,$\SS$\, of the mean curvature \,$\gH$\, of the cylinder, when \,$m=n^2$\, and \,$n\to \i$.
Conversely, it goes to \,$+\i$\, if \,$m=n$\, and \,$n\to\i$. \ep
\bpf We have to distinguish between the edges where two bases meet, %
and edges where two lateral sides of the triangles meet, respectively.
\par Concerning the \,$2m(n-1)$\, edges \,$e$\, where two bases meet, their length is \,$b$\, and all of them have the same dihedral angle %
\beq\label{teta1}\t_e=2\arctan\Bigl({d\over H/n}\Bigr)\,, \eeq
hence the mean curvature at each edge is \,$b\cdot\t_e/2$.
Therefore, the total contribution of the mean curvature at the first kind of edges is
$$ F_1(m,n):=2m(n-1)\cdot R\sin\a_m\cdot \arctan\Bigl({n\over H}\,2R\sin^2(\a_m/2) \Bigr)\,. $$
Since the bases of the triangles are orthogonal to the direction of the first principal curvature \,$\gk_1=0$\, of the cylinder, one expects that when suitably passing to the limit one gets \,$F_1(m,n)\to 0$. Actually, when \,$n=m$\, one computes \,$F_1(n,n)\to 2\p^3R^2/H$\, as \,$n\to \i$. However, taking e.g. \,$m=n^2$\, one gets
$$ \lim_{n\to\i}F_1(n^2,n)=\lim_{n\to\i}2n^2(n-1)\cdot R\sin\Bigl({\p\over n^2}\Bigr)\cdot \arctan\Bigl({n\over H}\,2R\sin^2\Bigl({\p\over 2n^2}\Bigr) \Bigr)=0\,. $$
\par We now deal with the \,$2mn$\, edges \,$e$\, where two lateral sides of the triangles meet. Such edges are almost orthogonal to the direction of the second principal curvature \,$\gk_2=1/R$\, of the cylinder. Therefore, one expects that when suitably passing to the limit this time one gets
\,$\p\,H$, i.e., the integral of the mean curvature of the cylinder. Any such edge has length equal to the lateral edge of the congruent triangles, whence
\,$|e|=\sqrt{h^2+(b/2)^2}$. Moreover, by the symmetry it turns out that all of them have the same dihedral angle \,$\t_e$, which will be computed by means of the formula
\beq\label{theta}  \t_e=\arcsin|\gn_1\wdg\gn_2| \eeq
where \,$\gn_1,\gn_2\in\Sph^2$\, are the outward unit normals of the two triangles meeting at the edge. We thus e.g. consider the first isosceles triangle with vertexes \,$A=(R,0,0)$, \,$B=(R\cos(2\a_m),R\sin(2\a_m),0)$, $C=(R\cos\a_m,R\sin\a_m,H/n)$. We compute
$$ \vect{AB}=R(-2\sin^2\a_m,2\sin\a_m\cos\a_m,0)=2R\sin\a_m\,(-\sin\a_m,\cos\a_m,0)\,,$$
$$|\vect{AB}|=b\,. $$
The middle point of the basis \,$AB$\, is \,$H=(R\cos^2\a_m,R\cos\a_m\sin\a_m,0)$, which gives
$$ \vect{HC}=(R\cos\a_m(1-\cos\a_m),R\sin\a_m(1-\cos\a_m),H/n)\,,\quad |\vect{HC}|=h\,. $$
Taking the wedge product of the orthogonal unit vectors \,$v_1:=\vect{AB}/|\vect{AB}|$\, and \,$v_2:=\vect{HC}/|\vect{HC}|$,
we get
$$ \gn_1:=v_1\wdg v_2={1\over h}\,\Bigl({H\over n}\cos\a_m,{H\over n}\sin\a_m,-R(1-\cos\a_m) \Bigr)\,,\quad |\gn_1|=1\,. $$
The adjacent triangle has vertexes \,$C$, $D$, and \,$A$, where $$D=(R\cos\a_m,-R\sin\a_m,H/n)\,,$$ whence
$$ \vect{DC}=(0,2R\sin\a_m,0)\,,\quad |\vect{DC}|=b\,. $$
The middle point of the basis \,$DC$\, is \,$K=(R\cos\a_m,0,H/n)$, which gives
$$ \vect{AK}=(-R(1-\cos\a_m),0,H/n)\,,\quad |\vect{AK}|=h\,. $$
Therefore, the wedge product of the orthogonal unit vectors \,$v_3:=\vect{DC}/|\vect{DC}|$\, and \,$v_4:=\vect{AK}/|\vect{AK}|$\, gives
$$ \gn_2:=v_3\wdg v_4={1\over h}\,\Bigl({H\over n},0,R(1-\cos\a_m) \Bigr)\,,\quad |\gn_2|=1\,. $$
\par Now, we have
$$ \gn_1\wdg\gn_2= {H\over n}\,{1\over h^2}\,{\bf v}\,,\quad {\bf v}:=\Bigl(R\,(1-\cos\a_m)\,\sin\a_m,-R\,\sin^2\a_m, \sin\a_m \Bigr)$$
where we compute \,$|{\bf v}|=h\,\sin\a_m$, so that we get
\beq\label{teta2}
 \t_e=\arcsin|\gn_1\wdg\gn_2|=\arcsin\Bigl({H\over n}\,{1\over h^2}\,\sin\a_m \Bigr)\,. \eeq
Therefore, the contribution to the mean curvature of \,$P_{m,n}$\, given by the \,$2mn$\, edges \,$e$\, where two lateral sides of the triangles meet is
\beq\label{F2} F_2(m,n)=2mn\cdot\sqrt{h^2+(b/2)^2}\cdot{\t_e \over 2} \eeq
where, we recall,
$$ b=2R\sin\a_m\,,\quad h=((H/n)^2+d^2)^{1/2}\,,\quad d=2R\sin^2(\a_m/2)\,,\quad \a_m={\p\over m}\,. $$
\par If \,$m=n$, it is readily checked that \,$F_2(n,n)\to +\i$\, as \,$n\to \i$.
Taking instead \,$m=n^2$, we have
\,$ 2mn\cdot\sqrt{h^2+(b/2)^2}\sim 2\,H\,n^2 $, whereas \,$\t_e\sim \a_{n^2}=\p/n^2$, whence
$$ \lim_{n\to \i}F_2(n^2,n)=\p\,H\,. $$
Since the integral of the mean curvature of the cylinder is $\int_\SS \gH\,d\Ha^2=\p\,H$, the proof is complete.
\epf
\begin{remark}\label{Rcurv} If we consider the curvature force term \,$2\sin(\t_e/2)$\, instead of the angle \,$\t_e$, it is readily seen that the above computation yields to the same conclusions.
We also notice that as for the area, and of course for the Gauss curvature, the lower semicontinuity property holds. Namely, letting \,$m,n\to\i$, the lower limit of the terms \,$F_1(m,n)+F_2(m,n)$\, is always greater that \,$\p\,H$. In fact, in the formula \eqref{F2} we have:
\,$2mn\sqrt{h^2+(b/2)^2}\geq 2mn h\geq 2mH$, whereas \,$\t_e\geq \arcsin[(1+(2R\,n/H)^2\sin^4\a_m )^{-1/2}\sin\a_m]$, which tends to
\,$\a_m$\, as \,$n\to\i$, and \,$mH\,\a_m=\p\,H$.
\end{remark}
\par\noindent{\large\sc The curvature energy of smooth surfaces.}
Following \cite{AST}, see also \cite{Anz}, for a smooth surface \,$\M$\, in \,$\gR^{3}$, all the information about the curvatures is contained in the graph
$$ \G\M:=\{(z,\n(z))\mid z\in \M\}$$
of the Gauss map \,$\n:\M\to\Sph^2\sb\gR^{3}$\, of the surface.
Since the tangent plane to \,$\G\M$\, at a point \,$(z,\n(z))$\, is
determined by the tangential derivatives of \,$\n(z)$\, at \,$z$, and hence by the second fundamental form to \,$\M$\, at \,$z$,
by the area formula it turns out that the area of the Gauss graph surface \,$\G\M$\, is linked to the principal curvatures
of \,$\M$\, by the relation:
\beq\label{ATS}
\Ha^2(\G\M)=\int_\M \Bigl( 1+
({\gk_1}^2+{\gk_2}^2)+(\gk_1\gk_2)^2 \Bigr)^{1/2}\,d\Ha^2 \eeq
where \,$\gk_1=\gk_1(z)$\, and \,$\gk_2=\gk_2(z)$\, are the {\em
principal curvatures} at \,$z\in\M$, and \,$\Ha^2$\, denotes the 2-dimensional Hausdorff measure. 
\par More precisely, the tangent 2-vector field \,$\tt:\M\to\Lambda^2
T\M\sb\LL^3\gR^3_z$\, is given in terms of the Hodge operator by
\,$\tt(z)=\ast \n(z)$. Denoting by
\,$\Phi:\M\to\gR^3_z\tim\gR^3_y$\, the graph map
\,$\Phi(z):=(z,\n(z))$, a continuous tangent 2-vector field
\,$\x:\G\M\to\bigwedge^2(\gR^3_z\tim\gR^3_y)$\, is given by
\,$\x(z,\n(z)):=\bigwedge^2 d\Phi_z(\tt(z))$.
Moreover, denoting by \,$\tt_1$\,
and \,$\tt_2$\, the {\em principal directions}, and considering
the obvious homomorphism \,$v\mapsto \wid v$\, from \,$\gR^3_z$\,
onto \,$\gR^3_y$, one has
\beq\label{ATS2}
\x(z,\n(z))=\tt_1\wdg\tt_2+ \Bigl(
\gk_2\tt_1\wdg\wid\tt_2 - \gk_1\tt_2\wdg\wid\tt_1\Bigr) +
\gk_1\gk_2\,\wid\tt_1\wdg\wid\tt_2 \eeq
and since \,$|\x|\geq
1$\, on \,$\G\M$, the normalized $2$-vector field
\,$\vect\z:=\x/|\x|$\, determines an orientation to \,$\G\M$.
The area formula gives
$$ \Ha^2(\G\M)=\int_\M J^\M_\Phi\,d\Ha^2$$
where \,$J^\M_\Phi(z)$\, is the tangential Jacobian to \,$\Phi$\, at \,$z$, see e.g. \cite[Sec.~2.11]{AFP}.
Using that \,$J^\M_\Phi(z)=|\x(z,\n(z))|$, formula \eqref{ATS} follows from \eqref{ATS2}.
%
%
\par Also, denoting by \,$\gH$\, and \,$\gK$\, the {\em mean curvature}
and {\em Gauss curvature},
$$ \gH:={1\over 2}(\gk_1+\gk_2)\,,\qquad \gK:=\gk_1\gk_2 $$
so that \,$\gk_{1,2}=\gH\pm\sqrt{\gH^2-\gK}$, one equivalently has
$$ ({J^\M_\Phi})^2=1+(2\gH)^2-2\gK+\gK^2=4\gH^2+(1-\gK)^2 $$
and hence
$$
\Ha^2(\G\M)=\int_\M
\sqrt{ 1+(4\gH^2-2\gK)+\gK^2 }\,d\Ha^2\,. $$
%
%
\par The {\em curvature functional} is defined in \cite{AST} by
$$ \Vert
\M\Vert:=\Ha^2(\M)+\int_\M\sqrt{\gk_1^2+\gk_2^2}\,d\Ha^2+\int_\M|\gk_1\gk_2|\,d\Ha^2
$$
i.e., equivalently,
\beq\label{vertM} \Vert
\M\Vert:=\int_\M\bigl( 1+\sqrt{4\gH^2-2\gK}+|\gK|\bigr)\,d\Ha^2\,.
\eeq
\par\noindent{\large\sc Non-parametric surfaces}.
We shall denote by \,$(e_1,e_2,e_3)$\, the canonical basis of
\,$\gR^3_z$. Also, for a function \,$v:Q\to\gR$, we shall always denote by
\,$\nabla v$\, the (approximate) gradient and by \,$\pa_i v$\, and
\,$\pa^2_{i,j} v$\, the first and second order (approximate) partial
derivatives, so that e.g. \,$ \pa_i v(x):=\nabla v(x)\bullet
e_i$\, for \,$i=1,2$.
\par Assume now \,$\M=\G_u$, where \,$\G_u:=\{z=(x,u(x))\mid x\in Q \} $ is the graph of
a smooth function \,$u:Q\to\gR$. The Gauss map is naturally identified at each point of the graph
by the outward unit normal
$$
\gnu(x):={1\over \sqrt{g_u}}\bigl(-\pa_1
u, -\pa_2 u, 1 \bigr)\,,\qquad g_u:=1+|\nabla u|^2,\qquad x\in Q $$
and hence the Gauss graph of the non-parametric smooth surface \,$\G_u$\, is
\,$\G\G_u=\{\Phi_u(x)\mid x\in Q \}$, where
\,$\Phi_u:Q\to \gR^3_z\tim\gR^3_y$\, is the smooth map
$$
\Phi_u(x)=\bigl(\vf_u(x),\gnu(x)\bigr)\,,$$
$$\vf_u(x):=(x,u(x))\,,\quad
\gnu(x)=(\gnu^1(x),\gnu^2(x),\gnu^3(x))\,. $$
The mean
curvature at \,$(x,u(x))$\, becomes
$$
\gH_u={1\over
2}\,{1\over {g_u}^{3/2}}\,\bigl(
(1+(\pa_1u)^2)\pa^2_{2,2}u+(1+(\pa_2u)^2)\pa^2_{1,1}u -2
\pa_1u\,\pa_2u\,\pa^2_{1,2}u\bigr) $$
and the Gauss curvature at \,$(x,u(x))$\,
$$
\gK_u={1\over
{g_u}^2}\,\Bigl( \pa^2_{1,1}u\,\pa^2_{2,2}u
-(\pa^2_{1,2}u)^2\Bigr)\,. $$
Therefore, by the area formula we can write the
area of the Gauss graph as
$$
\Ha^2(\G\G_u)= \int_{\G_u} \sqrt{ 1+(4\gH_u^2-2\gK_u)+\gK_u^2
}\,d\Ha^2=$$
$$=\int_Q\sqrt{ g_u}\sqrt{ 1+(4\gH_u^2-2\gK_u)+\gK_u^2
}\,dx\,. $$
\par The tangent space at each point in the Gauss graph \,$\G\G_u$\, is oriented by
the wedge product
$$ \x_u(x):=\pa_1\Phi_u(x)\wdg\pa_2\Phi_u(x)\,,\qquad x\in Q\,. $$
Let \,$(\e_1,\e_2,\e_3)$\, be the canonical basis in \,$\gR^3_y$, the ambient space of the unit normal \,$\n_u$. According to the number of \,$\e_j$-entries, we can
write as in \cite{AST} the {\em stratification}
$$ \x_u=\x_u^{(0)}+\x_u^{(1)}+\x_u^{(2)}$$
where, denoting by \,$|M|$\, the determinant of a $2\tim 2$ matrix
\,$M$, we compute:
\beq\label{strat}
\ba{l}\ds\x_u^{(0)}= e_1\wdg e_2+ \pa_2 u\,
e_1\wdg e_3-\pa_1u\, e_2\wdg e_3 \\
\ds \x_u^{(1)}=\sum_{j=1}^3\pa_2\gnu^j e_1\wdg\e_j -
\sum_{j=1}^3\pa_1\gnu^j e_2\wdg\e_j + \sum_{j=1}^3
\left\vert\ba{cc}
  \pa_1u & \pa_2u \\
  \pa_1\gnu^j & \pa_2\gnu^j
\ea\right\vert e_3\wdg\e_j \\
\ds\x_u^{(2)}= \left\vert\ba{cc}
  \pa_1\gnu^1 & \pa_2\gnu^1 \\
  \pa_1\gnu^2 & \pa_2\gnu^2
\ea\right\vert\e_1\wdg\e_2+
\left\vert\ba{cc}
  \pa_1\gnu^1 & \pa_2\gnu^1 \\
  \pa_1\gnu^3 & \pa_2\gnu^3
\ea\right\vert\e_1\wdg\e_3 + \left\vert\ba{cc}
  \pa_1\gnu^2 & \pa_2\gnu^2 \\
  \pa_1\gnu^3 & \pa_2\gnu^3
\ea\right\vert\e_2\wdg\e_3\,. \ea \eeq
By \eqref{ATS2}, we thus infer:
$$
|\x_u^{(0)}|^2=g_u\,,\quad
|\x_u^{(1)}|^2=g_u\,(4\gH_u^2-2\gK_u)\,,\quad
|\x_u^{(2)}|^2=g_u\,\gK_u^2 $$
and hence, by taking \,$\M=\G_u$\, in \eqref{vertM}, again by the area formula we can equivalently write the curvature functional \,$\Vert\G_u\Vert$\,
of a smooth non-parametric surface as
\beq\label{FuQ}
\F(u,Q):=A(u,Q)+\F_1(u,Q)+\F_2(u,Q) \eeq
where we have set
\beq\label{Eu}
A(u,Q):=\int_Q\sqrt{g_u}\,dx\,,$$
$$ \F_1(u,Q):=\int_Q \sqrt{g_u}\,\sqrt{4\gH_u^2-2\gK_u}\,dx\,,$$
$$ \F_2(u,Q):=\int_Q \sqrt{g_u}\,|\gK_u|\,dx\,. \eeq
We thus get:
%
$$
A(u,Q)=\int_Q |\x_u^{(0)}|\,dx\,,\quad\F_1(u,Q)=\int_Q |\x_u^{(1)}|\,dx\,,$$
$$
\F_2(u,Q)=\int_Q
|\x_u^{(2)}|\,dx\,,\quad \Ha^2(\G\G_u)=\int_Q |\x_u|\,dx
$$
where \,$|\x_u|^2=|\x_u^{(0)}|^2+|\x_u^{(1)}|^2+|\x_u^{(2)}|^2$\,
and more explicitly, by \eqref{strat},
\beq\label{xu} \ba{l} |\x_u^{(0)}|^2  =  g_u=1+|\nabla u|^2 \\
\ds |\x_u^{(1)}|^2  =  g_u\,(4\gH_u^2-2\gK_u) =
|\nabla\gnu|^2+\sum_{j=1}^3\bigl(\pa_1u\,\pa_2\gnu^j-\pa_2u\,\pa_1\gnu^j\bigr)^2
\\
\ds |\x_u^{(2)}|^2  =  g_u\,\gK_u^2 = \sum_{1\leq j_1<j_2\leq 3}
\bigl(
\pa_1\gnu^{j_1}\,\pa_2\gnu^{j_2}-\pa_2\gnu^{j_1}\,\pa_1\gnu^{j_2}\bigr)^2\,.
\ea \eeq
\begin{remark}\label{Rint} The three integrals \,$\int_Q|\x_u^{(i)}|\,dx$, for \,$i=0,1,2$, may be seen as the smooth counterpart of the energy terms \,$A(P)$, \,$\EE_\gH(P)$, and \,$\EE_\gK(P)$, respectively, for
polyhedral surfaces \,$P$\, inscribed in the graph of a continuous function as above. In fact, the first term, \,$A(u,Q)$, is equal to the area of the smooth graph surface \,$\G_u$. The second term, \,$\F_1(u,Q)$, depends on both the mean and Gauss curvature of \,$\G_u$, and actually \,$\F_1(u,Q)\geq\int_Q|\nabla\n_u|\,dx$, hence it provides an upper bound to the total variation of the smooth outward unit normal \,$x\mapsto\n_u(x)$.
The third term, \,$\F_2(u,Q)$, only depends on the Gauss curvature of \,$\G_u$, and by the area formula it agrees with the mapping area (i.e., counting the multiplicity) in \,$\Sph^2$\, of the outward unit normal. Finally, since \,$|\x_u|\leq (|\x_u^{(0)}|+|\x_u^{(1)}|+|\x_u^{(2)}| )\leq\sqrt{3}\,|\x_u|$, we have:
$$ {1\over\sqrt 3}\,\Ha^2(\G\G_u)\leq \F(u,Q):=A(u,Q)+\F_1(u,Q)+\F_2(u,Q)\leq \Ha^2(\G\G_u)$$
where, we recall, \,$\Ha^2(\G\G_u)$\, is the area of the Gauss graph \,$\G\G_u$\, of the graph surface \,$\G_u$.
\end{remark}
\section{Smoothing out a polyhedral chain}\label{Sec:smoothing}
In this section we analyze the curvature energy of smooth approximations of a polyhedral surface.
\par We thus assume that \,$v$\, is a Lipschitz function on \,$Q=[0,1]^2$\, which is affine on each triangle of a finite triangulation \,$\D$\, of the domain \,$Q$.
Then the graph of \,$v$\, is a triangulated polyhedral surface \,$P$.
Moreover, \,$\nabla v\in L^\i(Q)$\, and the unit normal \,$\n_v:Q\to\Sph^2$\, is a $BV$-function whose weak derivative \,$D\n_v$\, is a finite
vector-valued measure concentrated on the edges of the triangulation. By means of a convolution argument, we shall prove the following:
\bp\label{Pappr} There exists a sequence of smooth functions \,$u_h:Q\to\gR$\, such that \,$u_h$\, converges to \,$v$\, strongly in \,$W^{1,1}$, the unit normals \,$\nu_{u_h}$\, converge to \,$\n_v$\, strongly in the \,$\BV$-sense and finally
$$ \lim_{h\to\i}A(u_h,Q)=A(P)\,,\qquad \sup_h\F_1(u_h,Q)\leq {\p\over 2}\cdot \EE_\gH(P)\,. $$
More precisely, in terms of the energy \eqref{EPH2} we get:
$$\lim_{h\to\i}\F_1(u_h,Q)= \wid\EE_\gH(P)\,. $$ \ep
\par\noindent
{\large\sc Estimates by area in the Gauss sphere.}
In the proof of Proposition~\ref{Pappr}, in general it cannot be obtained a bound of the type
$$ \sup_h\F_2(u_h,Q)\leq C\cdot \EE_\gK(P)\,. $$
%
%
\bex\label{Ecylinder} Assume e.g. that \,$P$\, is (up to a rotation, so that \,$P$\, is the graph of a Lipschitz-continuous function) a piece, say \,$\wid P_{m,n}$, of the polyhedral surface \,$P_{m,n}$\,
obtained in the Schwarz--Peano example, in correspondence to a cylinder of radius \,$R$\, and height \,$H$.
In this example, at any vertex \,$V$\, in \,$P=\wid P_{m,n}$\, we know that \,$\gK_P(V)=0$, whence \,$\EE_\gK(P)=0$.
On the other hand, when looking at the smoothing argument, in a small neighborhood of each one of the six edges meeting at \,$P$, the outward unit normal of a smooth approximating function has to cover an arc in the Gauss sphere \,$\Sph^2$\, connecting the points given by the values of the outward unit normal to the two triangles of \,$P$\, meeting at the edge. The length of this arc is of the order of the dihedral angle \,$\t_e$\, at the edge \,$e$,
which is given by \eqref{teta1} for two edges, and by \eqref{teta2} for the other four ones.
%
On the other hand, see Remark~\ref{Rint}, the integral \,$\F_2(u_h,Q)$, that only depends on the Gauss curvature of the graph of \,$u_h$, is concentrated near the points of the 0-skeleton of the triangulation \,$\D$\, of \,$Q$\, that corresponds by projection to the triangles of \,$P$. Since the integral \,$\F_2(u_h,Q)$\, agrees with the mapping area of the outward unit normal of \,$u_h$, at each vertex it gives a contribution equal to the area (with multiplicity) of the spherical shell in \,$\Sph^2$\, enclosed by the ordered join of the six arcs previously described, which is a positive quantity, depending on \,$R,H,m$, and \,$n$. Proposition~\ref{Pestimatearea} below clarifies the situation, yielding to an upper bound of the area (with multiplicity) of the spherical shell in terms of the sum of the angles of the tiles concurring in the vertex.
\eex
\begin{remark}\label{Rvertexes} More generally, recalling \eqref{strat} and \eqref{Eu}, the area formula yields that for smooth functions, the energy \,$\EE_2(u,Q)$\, is equal to the mapping area of the unit normal \,$\n_u$\, in the sphere \,$\Sph^2$. If \,$\{u_h\}$\, is the smooth approximating sequence from Proposition~\ref{Pappr}, it turns out that the energy density
of the integral \,$\EE_2(u_h,Q)$\, is concentrated near the projection points in \,$Q$\, of the interior vertexes of the polyhedral surface \,$P$,
and around any such point the energy contribution is bounded (up to an absolute multiplicative constant) by the area of the {\em geodesical envelope}
of the unit normals of the triangular tiles of \,$P$\, concurring in the vertex. \end{remark}
\par\noindent
{\large\sc A rough estimate.} We shall prove the following area estimate:
\bp\label{Pestimatearea}
Let $V$ be a vertex of a polyhedral graph $P$. Let $N_0,\ldots,N_{k-1}$ be the unit normals (with positive $z$-component) of the tiles $\alpha_0,\ldots,\alpha_{k-1}$ concurring in $V$. Then the area of the geodesical envelope $G$ of $N_0,\ldots,N_{k-1}$ in the Gauss half-sphere satisfies
\begin{equation}\label{est}
A(G)\ \leq\ 2\pi\sum_{j\in\mathbb Z_k}\theta_j\,,
\end{equation}
where $\theta_j$ is the angle in $V$ of the tile $\alpha_j$.
\ep
\par As a consequence, if we define the Gauss energy of the polyhedral surface by
\beq\label{EPK2} \wid\EE_\gK(P):=\sum_{V\in P} |\wid\gK_{P}(V)|\,,\qquad |\wid\gK_{P}(V)|:=\sum_{i}\theta_i \eeq
where the summation is taken on all the vertexes of \,$P$, and \,$\t_i$\, is the angle of the $i^{th}$-triangle of \,$P$\,
meeting at \,$V$, we readily extend Proposition~\ref{Pappr} as follows:
\bc\label{Cappr} In Proposition~$\ref{Pappr}$, we also have:
$$\sup_h\F_2(u_h,Q)\leq C\cdot 2\p\cdot \wid\EE_\gK(P) $$
where \,$C>0$\, is an absolute constant, not depending on \,$v$, and \,$\wid\EE_\gK(P)$\, is given by \eqref{EPK2}. \ec
\par Recalling that the Gauss curvature at a vertex \,$V$\, is equal to the angle defect, i.e., \,$\gK_{P}(V):=2\pi-\sum_i\t_i $, we give the following
\bdf\label{Dvert} {\em We say
that \,$V$\, is an elliptic, parabolic, or hyperbolic vertex of \,$P$\,
if the angle defect is positive, zero, or negative, respectively.} \edf
\bex\label{EGauss} If \,$P$\, is the cylindrical surface from Example~\ref{Ecylinder}, each vertex is a parabolic one, the polyhedral surface being developable.
We thus have \,$\EE_{\gK}(P)=0$, but \,$\wid\EE_{\gK}(P)=N\cdot 2\p$, where \,$N$\, is the number (depending on \,$n$\, and \,$m$\,) of vertexes in \,$P$.
As a consequence, by Corollary~\ref{Cappr}, the approximating sequence from Proposition~\ref{Pappr} satisfies the energy bound
\,$\sup_h\F_2(u_h,Q)\leq C\cdot 4\p^2\,N$, depending on the number of vertexes, and hence it has nothing to do with the energy \,$\EE_\gK(P)$,
which is equal to zero.
\par Of course, a similar drawback occurs in presence of hyperbolic vertexes. \eex
%
\par\noindent{\large\sc Elliptic vertexes.} On the other hand, if $V$ is an elliptic vertex of a polyhedral graph $P$,
we can thus refine the estimate in Proposition~\ref{Pestimatearea}, showing that the Gauss curvature can be calculated in terms of a suitable area in the
Gauss sphere:
\bp\label{area=curvature}
Let $V$ be an elliptic vertex of a polyhedral graph $P$, i.e., \,$\gK_P(V)>0$. 
Let $N_0,\ldots,N_{k-1}$ be the unit normals with positive $z$-component of the tiles $\alpha_0,\ldots,\alpha_{k-1}$ concurring in $V$. Then the area of the geodesical envelope $G$ of $N_0,\ldots,N_{k-1}$ in the Gauss half-sphere equals the Gauss curvature concentrated in $V$:
\begin{equation}\label{area-curv}
A(G)\ =\ \gK_P(V)\ =\ 2\pi-\sum_{j=0}^{k-1}\theta_j\,,
\end{equation}
where $\theta_j$ is the angle in $V$ of the tile $\alpha_j$.
\ep
\par As a consequence, if all the vertexes of the polyhedral graph are of elliptic type, using Remark~\ref{Rvertexes} and Proposition~\ref{area=curvature},
we readily extend Proposition~\ref{Pappr} as follows:
\bc\label{Capprconv} In Proposition~$\ref{Pappr}$, assume that each vertex of the polyhedral graph \,$P$\, is of elliptic type.
Then we also have:
$$\sup_h\F_2(u_h,Q)\leq C\cdot\EE_\gK(P) $$
where \,$C>0$\, is an absolute constant, not depending on \,$v$, and \,$\EE_\gK(P)$\, is the Gauss curvature energy defined in \eqref{EPHK}.
\ec
\par\noindent{\large\sc Proofs.} We conclude this section by proving Propositions~\ref{Pappr}, \ref{Pestimatearea}, and \ref{area=curvature}.
\smallskip\par\noindent
\bpff{\sc of Proposition~\ref{Pappr}:} Let \,$v$\, be a Lipschitz-continuous function on \,$Q=[0,1]^2$\, which is affine on each triangle of a finite triangulation \,$\D$\, of the domain \,$Q$. Assume first for simplicity that \,$v$\, is constant in a small neighborhood of the boundary \,$\pa Q$\, of the square domain, and extend \,$v$\, in a constant way to \,$\gR^2\sm Q$. Let \,$\r:\gR^2\to\gR$\, be a smooth symmetric mollifier with support contained in the unit ball centered at the origin, and denote \,$\r_h(x)=h^{2}\r(h\,x)$\, for \,$h\in\Nat^+$.
Define \,$u_h:Q\to\gR^2$\, by
$$ u_h(x):=(\r_h\ast v)(x)=\int_{\gR^2}\r_h(y)\,v(x-y)\,d y\,,\qquad x\in Q\,. $$
\par Since \,$v$\, is differentiable a.e. on \,$\gR^2$, with approximate gradient \,$\nabla v\in L^\i$, for \,$i=1,2$\, we have \,$\pa_i u_h=\r_h\ast \pa_i v$\, and the sequence \,$u_h$\, strongly converges to \,$v$\, in \,$W^{1,\i}(Q)$, by dominated convergence. In particular,
$$ \lim_{h\to\i}\ds\int_Q\sqrt{1+|\nabla u_h|^2}\,dx= \int_Q\sqrt{1+|\nabla
v|^2}\,dx\,. $$
Moreover, using that \,$\Vert \nabla u_h\Vert_\i\leq \Vert v\Vert_\i$, it turns out that the sequence \,$\n_{u_h}$\, of unit normals converges to the unit normal \,$\n_v$\, {\em strongly in the $\BV$-sense}, i.e.,
\,$\n_{u_h}\to \n_v$\, strongly in
\,$L^1(Q,\gR^3)$\, and
$$ \lim_{h\to \i}\int_Q|\nabla \n_{u_h}|\,dx= |D\n_v|(Q) \,. $$
\par Denoting by \,$|M|$\, the determinant of a $2\tim 2$ real matrix \,$M$, consider now for each \,$j=1,2,3$\, the functions
\,$\m^j_{u_h}:Q\to \gR$\,
\beq\label{mhj}  \m^j_{u_h}(x):=\left\vert\ba{cc}
  \pa_1u_h(x) & \pa_2u_h(x) \\
  \pa_1\n_{u_h}^j(x) & \pa_2\n_{u_h}^j(x)
\ea\right\vert\,. \eeq
Since the sequence \,$\{|\nabla u_h|\}$\, is equibounded, whereas \,$\{\nabla\n_{u_h}^j\}$\, converges to \,$\nabla\n_v^j$\, strongly in $L^1$,
it turns out that the sequence \,$\{\m^j_{u_h}\}$\, is equibounded in \,$L^1(Q)$.
\par More precisely, the energy contribution of the integral of the functions \,$\m^j_{u_h}$\, concentrates at the interior edges \,$\wid e$\, of the 1-skeleton of the triangulation \,$\D$, and we claim that around any such edge the integral of
the energy density
\beq\label{mh} \m_{u_h}(x):=\bigl({|\nabla\n_{u_h}(x)|^2+{\m^1_{u_h}(x)}^2+{\m^2_{u_h}(x)}^2+{\m^3_{u_h}(x)}^2}\bigr)^{1/2}\,,\qquad x\in Q \eeq
is bounded (up to an absolute multiplicative constant) by the mean curvature
\,$|\gH_P(e)|$\, of \,$P$\, at the edge \,$e$\, that projects onto \,$\wid e$. Actually, we shall see that it converges to the product of the length \,$\calL(e)$\,
times the dihedral angle \,$\t_e$\, of the edge.
\par Recalling \eqref{xu}, on account of the first definitions from \eqref{EPHK} and \eqref{EPH2} we definitely obtain:
$$ \int_Q |\x_{{u_h}}^{(0)}|\,dx=\int_Q\sqrt{g_{u_h}}\,dx\to \int_Q\sqrt{1+|\nabla v|^2}\,dx=A(P)$$
$$ \int_Q |\x_{{u_h}}^{(1)}|\,dx=\int_Q \m_{u_h}(x) \,dx\to \wid\EE_\gH(P)\,,\qquad
 \sup_h\int_Q |\x_{{u_h}}^{(1)}|\,dx\leq {\pi\over 2}\cdot \EE_\gH(P)<\i\,. $$
\par Since the argument is local, the claim can be checked by considering (without loss of generality) the case when \,$\wid e$\, is parallel to the direction $e_2$, i.e.,
the wedge product of the unit normals \,$\gn_1$\, and \,$\gn_2$\, of the two triangles of \,$P$\, that meet at the edge \,$e$\, is a vector of the type \,$(0,\l,\m)$, where \,$\l\neq 0$.
In this case, inside the two triangles we must have \,$\nabla v\equiv (a,c)$\, and \,$\nabla v\equiv (b,c)$, respectively, for some real constants \,$a,b,c$.
By using the formula \eqref{theta}, it turns out that the dihedral angle at the edge \,$e$\, is
$$ \t_e=\arcsin\Bigl( {|b-a|\sqrt{1+c^2}\over \sqrt {1+c^2+b^2}\cdot\sqrt {1+c^2+a^2}} \Bigr)\,.  $$
Since moreover \,$\t_e$\, agrees with the angle between the two planar vectors \,$(\sqrt{1+c^2},a)$\, and \,$(\sqrt{1+c^2},b)$, we equivalently have:
\beq\label{arctan}  \t_e=\Bigl\vert\arctan\Bigl( {b\over \sqrt {1+c^2}}\Bigr)-\arctan\Bigl( {a\over \sqrt {1+c^2}}\Bigr)\Bigr\vert\,.  \eeq
\par For \,$h$\, large, we denote by \,$I_h(\wid e)$\, the open set given by the points in \,$Q$\, whose distance from the edge \,$\wid e$\, is smaller than \,$1/h$\,
and whose distance from the vertexes of the edge \,$\wid e$\, is greater than \,$1/h$. In \,$I_h(\wid e)$, it turns out that
the second derivative \,$\pa_2 u_h \equiv c$\, whereas the first derivative \,$\pa_1 u_h$\, only depends on the first variable \,$x_1$, and actually it takes values in the segment with end points \,$a$\, and \,$b$. This yields that \,$\pa_2\n_{u_h}\equiv 0$\, and \,$\pa_1\n_{u_h}$\, only depends on the first variable \,$x_1$. As a consequence, we have \,$|\nabla \n_{u_h}|=|\pa_1\n_{u_h}|$\, and also \,$\m_{u_h}^j= -c\cdot\pa_1\n^j_{u_h}$, for \,$j=1,2,3$, which yields:
$$ \m_{u_h}(x)=\sqrt{1+c^2}\cdot |\nabla \n_{u_h}(x)|\,,\qquad |\nabla \n_{u_h}(x)|=|\pa_1\n_{u_h}(x)|\,,\quad \fa\,x\in I_h(\wid e)\,.$$
Furthermore, since we have
$$ \n_{u_h}(x)={1\over \sqrt{1+c^2+f_h^2(x)}}\,(-f_h(x),c,1)\,,\qquad f_h(x):=\pa_1 u_h(x) $$
we readily compute on \,$I_h(\wid e)$
$$ |\nabla\n_{u_h}|={\sqrt{1+c^2}\cdot |\pa_1 f_h| \over 1+c^2+{f_h}^2 }=\Bigl|\pa_1\arctan\Bigl({f_h\over \sqrt{1+c^2}} \Bigr) \Bigr|\,.$$
As a consequence, using that \,$f_h=\pa_1 u_h$\, is equal to \,$a$\, and \,$b$\, on the lateral sides of the set \,$I_h(\wid e)$,
and denoting by \,$L$\, the length of \,$\wid e$, on account of formula \eqref{arctan} we get the estimate:
$$
\int_{I_h(\wid e)}\m_{u_h}(x)\,dx= \sqrt{1+c^2} \int_{I_h(\wid e)}|\nabla\n_{u_h}|\,dx\leq \sqrt{1+c^2}\cdot L\cdot\t_\e+o(1/h)\,,
$$
where \,$o(1/h)\to 0$\, as \,$h\to\i$. Finally, observing that \,$\sqrt{1+c^2}\cdot L$\, is equal to the length \,$\calL(e)$\, of the edge \,$e$, we have obtained:
$$ \lim_{h\to\i}\int_{I_h(\wid e)}\m_{u_h}(x)\,dx=\calL(e)\cdot \t_e\leq {\p\over 2}\cdot |\gH_P(e)|\,. $$
Since the energy of \,$\m_{u_h}$\, concentrates near the edges \,$e$, the claim is proved and the proof is complete under the additional assumption that \,$v$\, is constant in a small neighborhood of the boundary \,$\pa Q$.
\par In general, one has to argue similarly as above, but this time using a procedure as in the density result by Anzellotti--Giaquinta \cite{AG},
i.e., by stepping down the size of the mollification when going
to the boundary of \,$\pa Q$,  compare e.g. Thm.~1 in \cite[Sec.~4.1.1]{GMSl}.
We omit any further detail. \epff
\bpff{\sc of Proposition~\ref{Pestimatearea}:}
Let the tiles $\alpha_0,\ldots,\alpha_{k-1}$ be numbered in order around the vertex $V$.

Let us fix the index $j\in\mathbb Z_k$. Let $\alpha_{j-1},\alpha_{j},\alpha_{j+1}$ be three consecutive tiles, $N_{j-1},N_{j},N_{j+1}$ their normal vectors and $\theta_j$ the angle of the tile $\alpha_j$ in $V$. Consider the geodesical triangle $T_j$ of vertexes $N_{j-1},N_{j},N_{j+1}$. If we prove that its area is bounded by $2\pi\theta_j$ we are done, since the geodesical triangles $T_j$, $j\in\mathbb Z_k$, cover all of $G$ hence
$$A(G)\ \leq\ \sum_{j\in\mathbb Z_k} A(T_j)\ \leq\ 2\pi\sum_{j\in\mathbb Z_k}\theta_j\,,$$
i.e.\ the thesis.
Let us prove that
\beq\label{equatore}
A(T_j)\ \leq\ 2\pi\theta_j\,.
\eeq
The normal vectors $N_{j-1}$ and $N_j$ determine the direction of their common edge $e_{j-1}$ (which is $D_{j-1}=N_{j-1}\wedge N_j$). The edge $e_j$, being in the tile $\alpha_j$ has a direction $D_j\perp N_j$. Moreover the distance on the Gauss sphere between the directions $D_{j-1}$ and $D_j$ is precisely the angle $\theta_j$ between the edges $e_{j-1}$ and $e_j$. Since the edge $e_j$ belongs also to the tile $\alpha_{j+1}$, one has $N_{j+1}\perp D_j$. Hence $N_{j}$ and $N_{j+1}$ both belong to the same maximum circle of directions perpendicular to $D_j$ and their distance on the Gauss semi-sphere is at most $\theta_j$. Considering the equator $E_j$ containing $N_j$ and $N_{j-1}$, the geodesical triangle $T_j$ is all contained in the strip between $E_j$ and one of its parallels at a distance $\theta_j$. Since the area of this sector is less than $2\pi\theta_j$, estimate \eqref{equatore} follows.
\epff
\bpff{\sc of Proposition~\ref{area=curvature}:} The hypothesis of ellipticity at the vertex $V$, means that the geodesical polygon of vertexes $N_0,\ldots,N_{k-1}$ in the Gauss half-sphere is a geodesically convex polygon coinciding with the geodesical envelope $G$ of  $N_0,\ldots,N_{k-1}$.

By elementary spherical geometry, the area $A(G)$ of such a geodesical polygon is given by
\begin{equation}\label{areapolygon}
A(G)\ =\ \sum_{j=0}^{k-1}\gamma_j-(k-2)\pi\,,
\end{equation}
where $\gamma_j$ is the angle between vertexes $N_{j-1},N_j,N_{j+1}$. If we prove that $\gamma_j+\theta_j=\pi$ for every $j$, equation \eqref{areapolygon} reduces to \eqref{area-curv} and the proposition is proved.

In order to compute $\gamma_j+\theta_j$ we interpret these angles as geodesical arcs on the Gauss sphere.

As already observed in the proof of Proposition~\ref{Pestimatearea}, the number $\theta_j$ is the distance between $D_{j-1}$ and $D_j$ on the Gauss sphere, which are two points on the equator relative to the pole $N_j$ and perpendicular to $N_{j-1}$ and to $N_{j+1}$, respectively.

On the other hand, the angle $\gamma_j$ in $N_j$ is the distance between the points $E_j$ and $E_{j+1}$ on the same equator relative to $N_j$, on the geodesic arc connecting $N_j$ to respectively $N_{j-1}$ and $N_{j+1}$.

Let $F_j$ be the point opposite to $E_j$ on the Gauss sphere. Our thesis is equivalent to the fact that the distance $E_{j+1}F_j$ is equal to $\theta_j$, i.e., to the distance $D_jD_{j+1}$. Adding to both arcs the arc $E_{j+1}D_j$, it is equivalent to show that the arcs $F_jD_j$ and $E_{j+1}D_{j+1}$ are congruent, which is implied in turn by the fact that given any point $P\neq N_j$, the points on the equator of $N_j$ of type $E$ (on the geodesic connecting $N_j$ and $P$) and of type $D$ (perpendicular to $N_j$ and $P$) lie clearly at a right angle (a distance of $\pi/2$), regardless of the choice of $P$: indeed $D$ is a pole relative to the equator through $N_j$ and $P$, which contains also $E$.
\epff

\section{The Gauss curvature in Schwarz-Peano example}\label{Sec:SPGauss}

In this section, we consider again the lateral surface of a cylinder $\Sigma$ and the polyhedral surfaces $P_{m,n}$ given by $2mn$ congruent isosceles triangles, as defined in Sec.~\ref{Sec:rel}.
As we have seen, as $m$ and $n$ go to infinity the area of the polyhedral surfaces may or may not go to the area of the cylinder, depending on the relative rates of the two parameters going to infinity. We will show that the equality given by Proposition~\ref{area=curvature} for elliptic vertexes of a polyhedral surface drastically fails in this case, where the vertexes are of parabolic type, see Definition~\ref{Dvert}.

As the computation is much more complicated in this case, we put $R=H=1$, thus the principal curvatures are $\gk_1=0$ and $\gk_2=1$ and the Gauss curvature is zero. Since the polyhedral surface is developable, the Gauss curvature is equal to zero at each vertex $V$. But it cannot be estimated properly by the area in the Gauss sphere of the geodesical envelope of the normals of triangles concurring in $V$.

Indeed, let us put ourselves near a vertex $V$ of the polyhedral surface $P_{m,n}$ and estimate the area in the Gauss sphere of the geodesical envelope of the six normals ($N_j$, $j=1,\ldots,6$) to the six triangles concurring in $V$. Calculating the normals $N_j$ and approximating all the trigonometric functions that appear (keeping in mind that we'll let $m$ and $n$ go to infinity) one gets:
$$ \ba{ll} \ds N_1\ =\ \frac{(0,-m,\pi n)}{\sqrt{\pi^2n^2+m^2}} & \ds N_4\ =\ \frac{(0,m,\pi n)}{\sqrt{\pi^2n^2+m^2}} \\
\ds N_2\ =\ \frac{(-\alpha m^5,-\alpha m^5,n)}{\sqrt{2\alpha^2m^{10}+n^2}} & \ds N_5\ =\ \frac{(\alpha m^5,-\alpha m^5,-n)}{\sqrt{2\alpha^2m^{10}+n^2}} \\
\ds N_3\ =\ \frac{(-\alpha m^5,-\alpha m^5,-n)}{\sqrt{2\alpha^2m^{10}+n^2}} & \ds N_6\ =\ \frac{(\alpha m^5,-\alpha m^5,n)}{\sqrt{2\alpha^2m^{10}+n^2}} \ea $$
where $\alpha$ is a positive constant, not depending on $n$ and $m$. The two diagonals $N_2N_5$ and $N_3N_6$ are congruent. The area of the geodesical hexagon con be roughly estimated by the product of the Euclidean distances $d_1=|N_2-N_5|=|N_3-N_6|$ and $d_2=|N_1-N_4|$. The distance $d_1$ lies between $\sqrt{2}$ and $2$, regardless of $m$ and $n$:
$$\sqrt{2}\ <\ d_1 \ =\ 2\cdot\frac{\sqrt{\alpha^2m^{10}+n^2}}{\sqrt{2\alpha^2m^{10}+n^2}} \ <\ 2\,,$$
while
$$d_2\ =\ \frac{2m}{\sqrt{\pi^2n^2+m^2}}\,.$$
Thus if $m=n$, one has \,$\ds d_2=\frac{2}{\sqrt{\pi^2+1}}$\, and hence the area of the hexagon is greater than a positive constant independent of $m$. The same holds if $m=n^k$, for $k>1$.

On the other hand, if $n=m^k$, $k>1$, then
$$d_2 \ \simeq\ 2m^{1-k}$$
which indeed goes to zero. But the number of vertexes in the triangulation is $nm=m^{k+1}$, so the total area in the whole polyhedral surface is of the order of $m^2$ and definitely it diverges, as \,$m\to\i$.

\adl\par\noindent{\large\sc A smarter triangulation.} By choosing a different triangulation of $\Sigma$, it turns out that area, mean curvature, and Gauss curvature behave as expected. Indeed by inscribing a prism $Q_n$ with base a regular $n$-agon in $\Sigma$ and then triangulating the lateral faces of the prism as we like, we have:
\begin{enumerate}
\item the lateral area of the prism $Q_n$ goes to the area of $\Sigma$ as $n\to\infty$;
\item the mean curvature concentrated in each side of the prism (elsewhere the mean curvature vanishes) is $H$ times half the curvature at a vertex of the regular $n$-agon, i.e.
$$\frac12\left(\pi-\frac{(n-2)\pi}n\right)\ =\ \frac\pi{n}\,,$$
hence the total mean curvature on $Q_n$ is equal to $\pi H$, i.e., to the total mean curvature of \,$\Sigma$;
\item the Gauss curvature is zero at each vertex, since the triangulation is developable;
\item in the Gauss sphere, the area of the geodesical envelope of the normals of triangles concurring to a vertex $V$ is zero, since either all normals coincide (if $V$ is inside a face of the prism) or there are only two different normals (if $V$ is on an edge of the prism) and hence their geodesical envelope is an arc of geodesic.
\end{enumerate}

Thus it is possible to approximate the cylinder $\Sigma$ with inscribed polyhedral surfaces $Q_n$ in such a way to have that area, mean curvature and Gauss curvature go to those of $\Sigma$ as $n\to\infty$, and that the equality of Proposition~\ref{area=curvature} holds.

The Schwarz-Peano example shows that in general,
this procedure has to be done in a smart way, depending on the geometry of the surface \,$\SS$, as in general not all triangulations work properly.
\section{BV and measure properties}\label{Sec:BV}
In this section we analyze the structure properties of continuous functions with finite relaxed energy.
\par
The celebrated theorem by L. Tonelli asserts that the membership of a continuous function \,$u:Q\to\gR$\, to the class \,$\BV(Q)$\, is equivalent
to the existence of a sequence of piecewise affine functions uniformly converging to \,$u$\, and whose graphs have equibounded area, compare below.
We shall see that a similar statement holds true with our notion (cf. \eqref{Area}) of relaxed area:
\bp\label{PTonelli} Let \,$u:Q\to\gR$\, be continuous. Then \,$u\in\BV(Q)$\, if and only if \,$\A(u,Q)<\i$. \ep
\par\noindent
{\large\sc $\BV$-property.} Now, if \,$u$\, is a continuous function with finite relaxed energy \,$\EE(u,Q)$, see \eqref{Erel}, then it has finite relaxed area \,$\A(u,Q)$, whence \,$u$\, is a function of bounded variation.
As a consequence, the outward unit normal \,$\n_u$\, is well defined a.e. on \,$Q$\, by \,$\n_u:=g_u^{-1/2}(-\pa_1 u,-\pa_2u,1)$,
where \,$g_u:=1+|\nabla u|^2$, but in term of the approximate partial derivatives of \,$u$.
In this section we shall prove the following:
\bt\label{TBV} Let \,$u$\, be a continuous function with finite relaxed energy \,$\EE(u,Q)$, see \eqref{Erel}. Then the outward unit normal \,$\n_u:Q\to\Sph^2$\, is a function of bounded variation,
\,$\n_u\in\BV(Q,\gR^3)$. \et
\par For this purpose, we first point out that Proposition~\ref{Pappr} and a diagonal argument yield:
\bc\label{Cdiag} Let \,$u:Q\to\gR$\, be a continuous function with finite relaxed energy \,$\EE(u,Q)$. Then there exists a
sequence of smooth functions \,$u_h:Q\to\gR$\, such that \,$u_h$\, converges to \,$u$\, strongly in \,$L^1(Q)$, and
$$ \sup_h\bigl(\A(u_h,Q)+\F_1(u_h,Q)\bigr)\leq C\cdot \EE(u,Q) $$
for some absolute constant \,$C>0$, not depending on \,$u$. \ec
\par On account of Corollary~\ref{Cdiag}, the $\BV$-property in Theorem~\ref{TBV} will be
obtained below through a
slicing argument, by exploiting analogous results from \cite{AcMu} for the total curvature of Cartesian curves.
\adl\par\noindent{\large\sc The mean curvature energy term.} As it is clear in the smooth case, the $\BV$-property of the unit normal \,$\n_u$\, does not guarantee a bound to the (relaxed) energy corresponding to the mean curvature. For this purpose, we introduce a suitable class of distributions that retain all the information.
\par The {\em distributional divergence} of an $L^1$-vector field \,$\s:Q\to\gR^2$\,
is well-defined by duality through the formula
$$ \lan\DIV \s,\vf\ran:= -\int_Q\s(x)\bullet \nabla\vf(x)\,dx\,,\qquad \vf\in C^\i_c(\Circle Q)$$
where \,$\bullet$\, denotes the scalar product in \,$\gR^2$.
\bex\label{Esju} If \,$u:Q\to\gR$\, is a continuous function with finite relaxed energy, by
Proposition~\ref{PTonelli} the approximate partial derivatives $\pa_iu$ are summable functions in $Q$. Moreover, by Theorem~\ref{TBV}
the unit normal \,$\n_u$\, is a function in $L^\i(Q,\gR^3)$. We can thus define the vector fields \,$\s^j_u\in L^1(Q,\gR^2)$\, through the formula
$$ \s^j_u:=(-\n_u^j\,\pa_2 u,\,\n_u^j\,\pa_1 u)\,,\quad j=1,2,3\,. $$ \eex
\par When \,$u$\, is smooth, say of class \,$C^2$, using that \,$\pa^2_{1,2}u=\pa^2_{2,1}u$\, and integrating by parts, we get
$$ \lan\DIV \s^j_u,\vf\ran=\int_Q (\pa_1 u\,\pa_2\n_u^j-\pa_2 u\,\pa_1\n_u^j)\,\vf(x)\,dx$$
for each \,$\vf\in C^\i_c(\Circle Q)$, and hence \,$\DIV\s^j_u$\, is an absolute continuous signed measure
$$ \DIV \s^j_u=\Div \s^j_u\,\calL^2\pri Q
$$
with density equal to the pointwise divergence of \,$\s^j_u$. Moreover, we have \,$\Div\s^j_u(x)=\m^j_u(x)$\, for each \,$x\in Q$\,
where, according to \eqref{mhj},
\beq\label{muj}
\m^j_u(x):=\left\vert\ba{cc}
  \pa_1u(x) & \pa_2u(x) \\
  \pa_1\n_{u}^j(x) & \pa_2\n_{u}^j(x)
\ea\right\vert\,. \eeq
\par\noindent
{\large\sc Polyhedral surfaces.} For each continuous function \,$u:Q\to\gR$\, with finite relaxed energy,
it is well defined the vector-valued distribution \,$\gm_u:=(\gm_u^1,\gm_u^2,\gm_u^3)$, where
$$ \gm_u^j:= (D\n_u^j,\DIV\s^j_u)\,,\qquad j=1,2,3\,. $$
If the graph of $u$ is a polyhedral surface $P$, the distribution \,$\gm_u$\,
is a finite measure. Moreover, in this case the total variation of \,$\gm_u$\, is equal to the energy term \,$\EE_\gH(P)$\, in \eqref{EPHK}. More precisely, we have:
\bp\label{Pmeanpoly} Let \,$u$\, be continuous on \,$Q$\, and affine on each triangle of a finite triangulation \,$\D$\, of the domain \,$Q$. Then, the distribution \,$\gm_u$\,
is a measure concentrated on the 1-skeleton of \,$\D$, and we actually have:
\beq\label{meanpoly} |\gm_u|(\Circle Q) =\EE_\gH(P) \eeq
where \,$P$\, is the polyhedral surface given by the graph of \,$u$. \ep
\par\noindent
{\large\sc Measure property.}  More generally, we shall prove the following:
\bt\label{Tmeas} Let \,$u$\, be a continuous function with finite relaxed energy \,$\EE(u,Q)$, see \eqref{Erel}. Then
for each \,$j=1,2,3$\, the distributional divergence \,$\DIV \s^j_u$\, is a finite measure, i.e.,
$$ \sup\{ \lan\DIV \s^j_u,\vf\ran \mid \vf\in C^\i_c(\Circle Q)\,,\,\, \Vert\vf\Vert_\i\leq 1 \} \leq C\cdot \EE(u,Q)<\i $$
for some absolute constant \,$C>0$, not depending on \,$u$. Moreover, the decomposition
\beq\label{dec}
 \DIV \s^j_u=\m^j_u\,\calL^2\pri Q+(\DIV \s^j_u)^s
\eeq
holds, where \,$\m^j_u$\, is the summable function defined \,$\calL^2$-a.e. on \,$Q$\, by \eqref{muj}, and \,$(\DIV \s^j_u)^s$\, is singular w.r.t. the Lebesgue measure.
\et
\par\noindent
{\large\sc The Gauss curvature energy term.} As to the Gauss curvature energy, we do not have an analogous to Proposition~\ref{Pmeanpoly}. In fact, for polyhedral surfaces, the Gauss curvature can be represented as a sum of Dirac masses concentrated at the 0-skeleton of the triangulation. However, such Dirac masses cannot be seen as derivatives of functions depending on $u$. This can be checked if one considers the boundary of the current \,$GG_u$\, associated to the Gauss graph of \,$u$, see \cite[Sec.~7]{Mu} for further details.
%
%
%
\par However, as a consequence of Corollary~\ref{Capprconv} we have:
\bp\label{PGauss} Let \,$u$\, be a continuous function with finite relaxed energy \,$\EE(u,Q)$, see \eqref{Erel}. Let
$$  \left(
      \begin{array}{cc}
        \pa_1\n_u^1 & \pa_2\n_u^1 \\
        \pa_1\n_u^2 & \pa_2\n_u^2 \\
        \pa_1\n_u^3 & \pa_2\n_u^3 \\
      \end{array}
    \right)
$$
be the matrix of the approximate partial derivatives of the unit normal.
If \,$u$\, is strictly convex (or strictly concave), then all the $2\tim 2$-minors of the above matrix
are in \,$L^1(Q)$. \ep
\begin{remark} We expect that the claim in Proposition~\ref{PGauss} holds true without assuming strict convexity.
However, we are not able to prove this fact, due to the drawbacks illustrated in Example~\ref{EGauss} and concerning parabolic and hyperbolic vertexes of a
polyhedral surface. \end{remark}
\par\noindent
{\large\sc Tonelli's theorem.}
In the classical definition by Tonelli,
letting \,$I=[0,1]$, and denoting by \,$V_1(x_1)$\, and \,$V_2(x_2)$\, the total variation in \,$I$\, of the functions \,$u(x_1,\cdot)$\, and  \,$u(\cdot,x_2)$, respectively,
a function \,$u:Q\to\gR$\, has bounded variation provided that both the functions \,$x_i\mapsto V_i(x_i)$\, are summable in \,$I$.
In a modern sense, since \,$u\in L^1(Q)$, an equivalent property is requiring that the distributional partial derivatives \,$D_i u$\, are measures of finite total variation.
\par Following e.g. \cite{Co}, in one implication of Tonelli's theorem, one assumes the existence of a sequence \,$\{P_h\}$\, of polyhedral surfaces given by the graph of functions
\,$v_h:Q\to\gR$\, such that the sequence \,$\{v_h\}$\, converges to \,$u$\,
uniformly on \,$Q$, and such that \,$\sup_h A(P_h)=C<\i$.
\par For any test function \,$\vf\in C^\i_c(\Circle Q)$\, and for \,$i=1,2$, one has
$$ \lan D_iu,\vf\ran:=-\lan u,\pa_i\vf\ran=-\lim_{h\to\i}\lan v_h,\pa_i\vf\ran $$
whereas for each \,$h$
$$ |\lan v_h,\pa_i\vf\ran|\leq A(P_h)\cdot\Vert\vf\Vert_\i\leq C\cdot \Vert\vf\Vert_\i$$
where the real constant \,$C>0$\, is a uniform bound to the areas \,$A(P_h)$\, of the approximating polyhedral surfaces, yielding to the required property \,$|D_iu|(Q)\leq C<\i$\,
for \,$i=1,2$, and hence that \,$u\in\BV(Q)$.
\par On the other hand, the converse implication in Tonelli's theorem reduces to the following statement: if a continuous function \,$u$\, belongs to the class \,$\BV(Q)$,
then there exists a sequence \,$\{v_h\}$\, of continuous functions \,$v_h:Q\to\gR$\, which are affine on each triangle \,$\DD$\, of a finite triangulation \,$\D_h$\,
of the square domain, such that \,$v_h\to u$\, uniformly on \,$Q$\, and \,$\sup_h A(v_h,Q)<\i$.
\par In order to prove the above statement, firstly, by means of a convolution argument with a symmetric mollifier, one defines a smooth sequence \,$\{u_h\}\sb C^\i(Q)$\, which converges to \,$u$\, uniformly on \,$Q$\, and such that the integrals \,$\int_Q|\nabla u_h|\,dx$\, converge to the total variation \,$|Du|(Q)$.
Secondly, by the smoothness of \,$u_h$, for each \,$h$\, one can easily find a polyhedral surface as above such that \,$\Vert v_h-u_h\Vert_\i\leq 2^{-h}$\, and
\,$A(v_h,Q)\leq C\cdot\A(u_h,Q)$\, for some absolute constant \,$C>0$, not depending on \,$h$. A diagonal argument yields the assertion.
\adl\par\noindent
{\large\sc Proofs.} We now give the proofs of Proposition~\ref{PTonelli}, Theorem~\ref{TBV}, Proposition~\ref{Pmeanpoly}, Theorem~\ref{Tmeas},
and Proposition~\ref{PGauss}.
\adl\par\noindent \bpff{\sc of Proposition~\ref{PTonelli}:} One implication follows by arguing as above. In fact, if the graph of \,$v_h$\, is an inscribed polyhedral surface \,$P_h=P_h(u,\D_h)$\, generated by the values $(x,u(x))$\, at the points \,$x$\,
in the 0-skeleton of a triangulation \,$\D_h$\, of the domain \,$Q$, condition \,$\mesh \D_h\to 0$\, yields that \,$v_h$\, uniformly converges to \,$u$.
As to the converse implication, if \,$\{v_h\}$\, is the sequence in Tonelli's theorem, we may and do assume that $\mesh \D_h\to 0$.
Letting \,$P_h=P_h(u,\D_h)$, by uniform convergence of \,$v_h\to u$, one infers (possibly passing to a subsequence) that
\,$\sup_h A(P_h)\leq C\cdot\sup_h A(v_h,Q)$\, for some absolute constant \,$C>0$. We omit any further detail. \epff
%
%
%
\bpff{\sc of Theorem~\ref{TBV}:} If \,$u:Q\to\gR$\, is continuous and with finite relaxed energy \,$\EE(u,Q)$, by Corollary~\ref{Cdiag} we may choose a smooth sequence \,$u_h:Q\to\gR$\, strongly converging to \,$u$\, in \,$L^1(Q)$\, and such that
\,$ \sup_h\bigl(\A(u_h,Q)+\F_1(u_h,Q)\bigr)\leq C<\i$. By Tonelli's theorem, we already know that \,$u\in\BV(Q)$, whence the outward unit normal \,$\n_u$\, is well defined \,$\calL^2$-a.e. in \,$Q$\, in terms of the approximate gradient of \,$u$, so that \,$\n_u\in L^1(Q,\Sph^2)$. Moreover,
by \eqref{Eu} we have for each \,$h$
$$
\F_1(u_h,Q)=\int_Q \sqrt{g_{u_h}}\,\sqrt{4\gH_{u_h}^2-2\gK_{u_h}}\,dx=\int_Q
|\x_{u_h}^{(1)}|\,dx $$
and hence by the second line in \eqref{strat} we get
$$ \sup_h\int_Q|\nabla \n_{u_h}|\,dx\leq\sup_h\F_1(u_h,Q)<\i\,. $$
As a consequence, by closure-compactness, see \cite{AFP}, possibly passing to a (not relabeled) subsequence, it turns out that the sequence \,$\{\n_{u_h}\}$\, weakly converges in the $\BV$-sense to some map \,$w\in\BV(Q,\Sph^2)$.
\par We now claim that for \,$i=1,2$, and possibly passing to a (not relabeled) subsequence, {\em the partial derivatives \,$\pa_i u_h$\, converge \,$\calL^2$-a.e. in \,$Q$\, to the approximate partial derivative \,$\pa_i u$}. This property implies that the corresponding sequence \,$\{\n_{u_h}\}$\, converges \,$\calL^2$-a.e. in \,$Q$\, to the outward unit normal \,$\n_u$. We thus have \,$w=\n_u$\, and hence \,$\n_u\in\BV(Q,\Sph^2)$. We recall, in fact, that by lower-semicontinuity of the total variation w.r.t. the weak \,$\BV$-convergence, one has
\,$|D\n_u|(Q)\leq\liminf_h\int_Q|\nabla \n_{u_h}|\,dx<\i$.
\par In order to prove the claim for e.g. \,$i=1$, letting \,$I=[0,1]$, for each \,$x_2\in I$\, and \,$t\in I$\, we shall denote \,$u^{x_2}_h(t):=u_h(t,x_2)$, and consider the smooth Cartesian curve \,$c^{x_2}_h(t):=(t,u^{x_2}_h(t))$.
We have:
\beq\label{int1}
\int_I\TC(c^{x_2}_h)\,dx_2\leq \A(u_h,Q)+\F_1(u_h,Q)\qquad \fa\,h\,. \eeq
\par In fact, using that \,$\dot c^{x_2}_h(t)=(1,\pa_1u_h(t,x_2))$\, we get
\,$g_{u_h}(t,x_2)\geq |\dot c^{x_2}_h(t)|^2$. Furthermore, recalling that the term \,$4\gH_{u_h}^2(t,x_2)-2\gK_{u_h}(t,x_2)$\, is equal to the sum of the square of the principal curvatures to the graph surface \,$\G_{u_h}$\, at the point \,$(t,x_2,u(t,x_2))$, such a quantity is greater than the square of the curvature
\,$\gk_{c^{x_2}_h}$\, of the curve \,$c^{x_2}_h$\, at the point \,$c^{x_2}_h(t)$, where
$$ \gk_{c^{x_2}_h}(t)={\pa^2_{1,1}u_h(t,x_2)\over (1+\pa_{1}u_h(t,x_2)^2)^{3/2}}\,,\quad t\in I\,.
$$
Since by the area formula
$$\TC(c^{x_2}_h)=\int_{c^{x_2}_h}|\gk_{c^{x_2}_h}|\,d\Ha^1 = \int_I|\dot c^{x_2}_h(t)|\,|\gk_{c^{x_2}_h}(t)|\,dt $$
we get
$$ \ba{rl} \ds \int_I\TC(c^{x_2}_h)\,dx_2= & \ds \int_I\Bigl(\int_I|\dot c^{x_2}_h(t)|\,|\gk_{c^{x_2}_h}(t)|\,dt\Bigr)dx_2 \\
\leq& \ds \int_Q \sqrt{g_{u_h}}\sqrt{1+4\gH_{u_h}^2-2\gK_{u_h}}\,dt\,dx_2 \ea$$
and hence inequality \eqref{int1} follows from the definitions in \eqref{Eu}.
\par By Fatou's lemma, we thus get
$$ \ba{rl}\ds\int_I\liminf_{h\to\i}\TC(c^{x_2}_h)\,dx_2 \leq & \ds\liminf_{h\to\i}\int_I\TC(c^{x_2}_h)\,dx_2 \\ \leq & \ds\sup_h \bigl(\A(u_h,Q)+\F_1(u_h,Q)\bigr)\leq C<\i\,. \ea $$
Moreover, following the notation from \cite{AcMu}, and letting \,$\ds\tt^{x_2}_{u_h}(t):={\dot c^{x_2}_{u_h}(t) \over |\dot c^{x_2}_{u_h}(t)|}$, we observe that
$$|\dot\tt^{x_2}_{u_h}(t)|=|\dot c^{x_2}_{u_h}(t)|\,|\gk_{c^{x_2}_{u_h}}(t)|\qquad \fa\, t\in I  $$
and hence we have:
$$ \int_I\liminf_{h\to\i}\int_I |\dot\tt^{x_2}_{u_h}(t)|\,dt\, dx_2\leq C<\i\,. $$
Therefore, by Rellich's theorem, we can find a (not relabeled) subsequence, not depending on \,$x_2\in I$, such that for \,$\calL^1$-a.e. \,$x_2\in I$\,
the sequence \,$\{\tt^{x_2}_{u_h} \}$\, converges weakly in the $\BV$-sense to some function \,$w^{x_2}:I\to\gR^4$. Now, arguing as in the proof of
\cite[Thm. 5.7]{AcMu}, we deduce that for \,$\calL^1$-a.e. $t\in I$
$$ w^{x_2}(t)={\dot c^{x_2}_{u}(t) \over |\dot c^{x_2}_{u}(t)|}\,, \qquad c^{x_2}_{u}(t):=(t,\pa_1u(t,x_2)) \,. $$
Since for \,$\calL^1$-a.e. \,$x_2\in I$\, the sequence \,$\{\tt^{x_2}_{u_h} \}$ converges \,$\calL^1$-a.e. in \,$I$\, to \,$w^{x_2}$, arguing as in \cite[Cor. 5.9]{AcMu} we conclude that the first partial derivative \,$\pa_1u_h(t,x_2)$\,
converges \,$\calL^2$-a.e. in \,$Q$\, to the approximate first partial derivative \,$\pa_1u_h(t,x_2)$, as required. \epff
\bpff{\sc of Proposition~\ref{Pmeanpoly}:} The first statement is trivial.
In order to check the equality \eqref{meanpoly}, choose an element \,$\wid e$\, of the 1-skeleton of \,$\D$, and assume (without loss of generality) that it is parallel to the direction $e_2$.
Therefore, inside the two triangles \,$\DD_1,\,\DD_2$\, of \,$\D$\, that meet at \,$\wid e$\, we have \,$\nabla u\equiv (a_1,c)$\, and \,$\nabla u\equiv (a_2,c)$,
respectively, so that for \,$i=1,2$\,
$$\n_{u\vert\Circle\DD_i}=\n_i:={1\over \sqrt{1+c^2+a_i^2}}\,(-a_i,c,1) \,. $$
If \,$e$\, is the edge in \,$P$\, that projects onto \,$\wid e$, we thus have \,$\calL(e)=|\wid e|\cdot\sqrt{1+c^2}$, whereas
$$ 2\,\sin\Bigl({\t_e \over 2}\Bigr)=\sqrt 2\,\sqrt{1-\cos\t_e}\,,$$
$$\cos\t_e=\n_1\bullet\n_2
= {1\over \sqrt{1+c^2+a_1^2}}\cdot{1\over \sqrt{1+c^2+a_2^2}}\,(1+c^2+a_1a_2)\,.$$
On the other hand, we have
$$(\s_u^1,\s_u^2,\s_u^3)_{u\vert\Circle\DD_i}={1\over \sqrt{1+c^2+a_i^2}}\,\bigl( -a_i\,(-c,a_i),-c\,(-c,a_i), (-c,a_i) \bigr) $$
so that on the line segment \,$\wid e$\, we compute:
$$ \ba{rl} |\gm^1_u|(\wid e)= & \ds \sqrt{1+c^2}\cdot\Bigl\vert{a_1\over \sqrt{1+c^2+a_1^2}}-{a_2\over \sqrt{1+c^2+a_2^2}}\Bigr\vert \\
|\gm^2_u|(\wid e)= & \ds \sqrt{1+c^2}\cdot\Bigl\vert{c\over \sqrt{1+c^2+a_1^2}}-{c\over \sqrt{1+c^2+a_2^2}}\Bigr\vert \\
|\gm^3_u|(\wid e)= & \ds \sqrt{1+c^2}\cdot\Bigl\vert{1\over \sqrt{1+c^2+a_1^2}}-{1\over \sqrt{1+c^2+a_2^2}}\Bigr\vert \ea $$
and definitely we get:
$$ |\gm_u|(\wid e)= \sqrt{1+c^2}\cdot\sqrt 2\,\sqrt{1-\cos\t_e}$$
which yields \eqref{meanpoly}, on account of \eqref{EPHK}. \epff
\bpff{\sc of Theorem~\ref{Tmeas}:} Let \,$\{u_h\}$\, be the smooth sequence given by Corollary~\ref{Cdiag}. Since
$$ \DIV \s^j_{u_h}=\Div \s^j_{u_h}\,\calL^2\pri Q\,,\quad\Div\s^j_{u_h}(x)=\m^j_{u_h}(x) \quad\fa\,x\in Q$$
where \,$\m^j_{u_h}(x)$\, is given by \eqref{mhj}, whereas by \eqref{Eu}, \eqref{xu}, and \eqref{mh} we estimate
$$ \int_Q |\m^j_{u_h}(x)|\,dx\leq \int_Q |\m_{u_h}(x)|\,dx=\int_Q |\x_{u_h}^{(1)}(x)|\,dx=\F_1(u_h,Q) $$
for \,$j=1,2,3$, possibly passing to a subsequence we deduce that the sequence of measures \,$\DIV \s^j_{u_h}$\, weakly converges to a signed measure \,$\gm^j$.
By lower-semicontinuity, we have
$$ |\gm^j|(Q)\leq \liminf_{h\to\i}|\DIV \s^j_{u_h}|(Q)=\liminf_{h\to\i}\int_Q |\m^j_{u_h}(x)|\,dx \leq C\cdot \EE(u,Q)<\i$$
and hence \,$\gm^j$\, has finite total variation.
We now claim that the following decomposition holds:
$$ \gm^j=\DIV\s^j_u+(\gm^j)^s $$
where the component \,$(\gm^j)^s$\, is singular with respect to the Lebesgue measure. In fact, as a consequence we also get
$$ |\DIV\s^j_u|(Q)\leq |\gm^j|(Q)\leq C\cdot \EE(u,Q)<\i\,.$$
\par In order to prove the claim, we recall that in the proof of Theorem~\ref{TBV} we have shown that, possibly passing to a (not relabeled) subsequence:
\ben \item the sequence of gradients \,$\nabla u_h$\, converge \,$\calL^2$-a.e. in \,$Q$\, to the approximate gradient \,$\nabla u$\,;
\item the sequence \,$u_h$\, weakly converges in the $\BV$-sense to \,$u\in \BV(Q)$\,;
\item
the sequence of unit normals \,$\{\n_{u_h}\}$\, weakly converges in the $\BV$-sense to the unit normal \,$\n_u\in\BV(Q,\Sph^2)$.
\een
\par By these properties, we have \,$\m^j_{u_h}\to\m^j_u$\, for \,$\calL^2$-a.e. \,$x\in Q$, hence by Fatou's Lemma
$$ \int_Q|\m^j_u(x)|\,dx\leq\liminf_{h\to\i}\int_Q|\m^j_{u_h}(x)|\,dx<\i $$
which yields that \,$\m^j_u\in L^1(Q)$\, for \,$j=1,2,3$. 
\par Moreover, for each \,$\vf\in C^\i_c(\Circle Q)$, possibly passing to a subsequence we deduce that \,$\s^j_{u_h}(x)\to \s^j_{u}(x)$\,
for \,$\calL^2$-a.e. \,$x\in Q$, whence:
$$  \ba{rl} \ds \lan\gm^j,\vf\ran= & \ds\lim\lm_{h\to\i}\lan \DIV\s^j_{u_h},\vf\ran=-\lim\lm_{h\to\i}\int_Q \s^j_{u_h}\bullet\nabla\vf\,dx \\
= & \ds -\int_Q \s^j_{u}\bullet\nabla\vf\,dx+\lan(\gm^j)^s,\vf\ran \ea $$
where the measure \,$(\gm^j)^s$\, is singular w.r.t. the Lebesgue measure \,$\calL^2\pri Q$.
The claim follows on account of the definition of distributional divergence \,$\DIV\s^j_{u}$. On the other hand, for each \,$h$
$$ -\int_Q \s^j_{u_h}\bullet\nabla\vf\,dx=\int_Q \m^j_{u_h}\cdot\vf\,dx $$
whereas, we recall, \,$\m^j_{u_h}\to\m^j_u$\, for \,$\calL^2$-a.e. \,$x\in Q$. This implies the decomposition \eqref{dec}, as required.
\epff
\bpff{\sc of Proposition~\ref{PGauss}:} If \,$u$\, is a strictly convex (or concave) function with finite relaxed energy, and \,$P$\, is a polyhedral graph inscribed in the graph of \,$u$,
it turns out that each vertex of \,$P$\, is of elliptic type.
Therefore, by Corollary~\ref{Capprconv}, and by a diagonal argument, we can find a sequence of smooth functions \,$u_h:Q\to\gR$\, uniformly converging to \,$u$\,
and such that \,$\sup_h\F_2(u_h,Q)\leq C\cdot\EE(u,Q)$, where \,$C>0$\, is an absolute constant.
The claim readily follows on account of \eqref{Eu} and \eqref{xu}, by lower-semicontinuity and by
the a.e. convergence of \,$\nabla\n_{u_h}$\, to the approximate gradient \,$\nabla\n_u$. Notice that the last property can be checked by means of the structure properties of the weak limit of the currents carried by the Gauss graph \,$\G\G_{u_h}$\, of the smooth approximating sequence, see \cite[Thm.~3.4]{Mu} and Sec.~\ref{Sec:rel}.
\epff
\adl\par\noindent{\bf Acknowledgements.}
{The research of D.M. was partially supported by PRIN 2010-2011 ``Calcolo delle Va\-ria\-zio\-ni'' and by the GNAMPA of INDAM.
The research of A.S. was partially supported by PRIN 2010-2011 ``Variet\`a reali e complesse:
geometria, to\-po\-lo\-gia e analisi ar\-mo\-ni\-ca'' and by the GNSAGA of INDAM.
We wish to thank the referee for his or her helpful remarks which helped to increase the readability of the paper.}
\end{document}